\documentclass[a4paper]{article}

\usepackage{amsfonts,amssymb} 
\usepackage[pdftex]{graphicx} % for \includegraphics
\usepackage[dvipsnames]{xcolor} % use predefined colors via name, e.g, "cyan"
\usepackage[normalem]{ulem} % for strike through
\usepackage{mathtools, mathabx} % extending amsmath
\usepackage{dsfont}
\usepackage[utf8]{inputenc}
\usepackage[T1]{fontenc} %The T1 font encoding is an 8-bit encoding and uses fonts that have 256 glyphs. So an 'ö' is an actual single glyph in the font. 
\usepackage[margin=1cm]{caption} % margin for captions
\usepackage[rightcaption,raggedright]{sidecap}% for side captions
\usepackage{framed} % framing boxes
\usepackage{float} % force table and figures with h, t, ..
\usepackage[section]{placeins} % forces floating objects (figures..) to remain within a section
\usepackage{algorithm} 
\usepackage{algpseudocode} 
\usepackage{hyperref}
\usepackage{chngcntr}
\counterwithin{figure}{section}
\usepackage{subcaption}
\usepackage{enumitem}

\usepackage{longtable}
\usepackage{booktabs}

\usepackage{color}%May be necessary if you want to color links
\definecolor{codegreen}{rgb}{0,0.6,0}

\usepackage{hyperref}
\hypersetup{
	colorlinks=true, %set true if you want colored links
	linktoc=all,     %set to all if you want both sections and subsections linked
	linkcolor=blue,  %choose some color if you want links to stand out
	citecolor=green,
	urlcolor=blue
}

\usepackage{amsthm}

\newtheorem{theorem}{Theorem}[section]

\newtheorem{definition}[theorem]{Definition}

\newboolean{is_comments}   
\setboolean{is_comments}{true}  

\newcommand{\holdAll}{\widehat{\Omega}}
\newcommand{\base}{\widetilde{\Omega}}

\newcommand{\ind}{\mathds{1}}
\newcommand{\smalltop}{{\!\top\!\!}}
\newcommand{\toPhys}[2]{\chi_{#1}(#2)}

\newcommand{\localStiffnessMatrix}[2]{A^{\#}_{{#1} {#2}}}
\newcommand{\localStiffness}[4]{A^{\#}_{{#1} {#2}} (\phi_{#3}, \phi_{#4})}
\newcommand{\termNonloc}[4]{A^{2}_{{#1} {#2}} ( \phi_{#3}, \phi_{#4})}%(k, k', j, j') }
\newcommand{\termLocal}[4]{A^{1}_{{#1} {#2}} ( \phi_{#3}, \phi_{#4})}%(k, k', j, j') }
\newcommand{\termNonlocPrime}[4]{A^{4}_{{#1} {#2}} ( \phi_{#3}, \phi_{#4})}%(k, k', j, j') }
\newcommand{\termLocalPrime}[4]{A^{3}_{{#1} {#2}} (  \phi_{#3}, \phi_{#4})}

\newtheorem{assumption}{Assumption}
\theoremstyle{definition}
\newtheorem{rem}[theorem]{Remark}

\newcommand{\InfinityNorm}{infinity norm}
\newcommand{\ApproximateCaps}{\emph{approxcaps} }
\newcommand{\NoCaps}{\emph{nocaps} }
\newcommand{\commonQuadRule}{\heartsuit}
\newcommand{\noZeroDivisionQuadRule}{\diamondsuit}
\newcommand{\duffyQuadRule}{\clubsuit}
\newcommand{\MK}[1]{\ifthenelse{\boolean{is_comments}}{\textcolor{blue}{[{#1}]}}{}  }
\newcommand{\CV}[1]{\ifthenelse{\boolean{is_comments}}{\textcolor{cyan}{[{CV:#1}]}}{}  }
\newcommand{\VS}[1]{\ifthenelse{\boolean{is_comments}}{\textcolor{red}{[{VS:#1}]}}{}  }
\newcommand{\doneCV}[1]{\ifthenelse{\boolean{is_comments}}{\textcolor{red}{[\sout {#1}]}}{}}
\newcommand{\doneMK}[1]{\ifthenelse{\boolean{is_comments}}{\textcolor{blue}{[\sout {#1}]}}{}}

\newcommand{\R}{\mathbb{R}}

\newcommand{\mbT}{\mathbb{T}}

\newcommand{\mcK}{\mathcal{K}}

\newcommand{\mcE}{\mathcal{E}}

\newcommand{\mcN}{\mathcal{N}}

\newcommand{\mcL}{\mathcal{L}}

\newcommand{\mcT}{\mathcal{T}}
\newcommand{\mcO}{\mathcal{O}}
\newcommand{\mcP}{\mathcal{P}}

\def \Cb{\mathbf{C}}

\def \ddim{d}

\def \fb{\mathbf{f}}

\def \gb{\mathbf{g}}

\def \idx{i}
\def \Jdx{J}
\def \jdx{j}
\def \Kdx{K}
\def \kdx{k}

\def \Mdx{M}

\def \ldx{\ell}
\def \ndim{n}

\def \pdx{p}
\def \Pdx{P}

\def \qdx{q}
\def \Qdx{Q}

\def \ub{\mathbf{u}}

\def \vb{\mathbf{v}}

\def \xb{\mathbf{x}}
\def \xbh{\widehat{\mathbf{x}}}

\def \yb{\mathbf{y}}
\def \ybh{\widehat{\mathbf{y}}}

\def \quadxp{d_{{\xbh}_\pdx}}
\def \quadyq{d_{{\ybh}_\qdx}}

\def \be{\begin{equation}}
\def \ee{\end{equation}}

\DeclareMathOperator{\spann}{span}

\DeclareMathOperator{\conv}{conv}

\DeclareMathOperator{\interior}{int}
\DeclareMathOperator{\Div}{div}
\title{\textbf{\texttt{nlfem}: A flexible 2d FEM Code for Nonlocal Convection--Diffusion and Mechanics}}

\author{
	Manuel Klar 
	\thanks{Universitaet Trier, D-54286 Trier, Germany; Email: klar@uni-trier.de, volker.schulz@uni-trier.de, vollmann@uni-trier.de} 
	\and 
	Christian Vollmann\footnotemark[1]
	\and
	Volker Schulz\footnotemark[1]}

\date{}

\let\oldproof\proof
\renewcommand{\proof}{\color{black}\oldproof}

\begin{document}
%%		
%%%%%%%%%%%%%%%%%%%%%%%%%%%%%%%%%%%%%%%%%%%%%%%%%%%%%%%%%%%%%%%%%%%%%%%%%%%%%%%%%
% HEADER
%%%%%%%%%%%%%%%%%%%%%%%%%%%%%%%%%%%%%%%%%%%%%%%%%%%%%%%%%%%%%%%%%%%%%%%%%%%%%%%%%
\maketitle
%ABSTRACT
\small
\textbf{Abstract.} In this work we present the mathematical foundation of an assembly code for finite element approximations of nonlocal models with compactly supported, weakly singular kernels. We demonstrate the code on a nonlocal diffusion model in various configurations and on a two-dimensional bond-based peridynamics model. Further examples can already be found in \cite{approxBall}. The code \texttt{nlfem} is published under the \texttt{MIT} License\footnote{For details see, e.g., \href{https://opensource.org/licenses/MIT}{https://opensource.org/licenses/MIT}.} and can be freely downloaded at \href{https://gitlab.uni-trier.de/pde-opt/nonlocal-models/nlfem}{https://gitlab.uni-trier.de/pde-opt/nonlocal-models/nlfem}.\\ 

%KEYWORDS
\normalsize

%%%%%%%%%%%%%%%%%%%%%%%%%%%%%%%%%%%%%%%%%%%%%%%%%%%%%%%%%%%%%%%%%%%%%%%%%
% CONTENT
%%%%%%%%%%%%%%%%%%%%%%%%%%%%%%%%%%%%%%%%%%%%%%%%%%%%%%%%%%%%%%%%%%%%%%%%%

% INTRODUCTION
\section{Introduction}\label{sec:introduction}
%1 why is this problem interesting
Over the last two decades nonlocal models attracted attention due to their capability of circumventing limitations of classical, local models based on differential equations. While the theory and application of differential equations is well established, nonlocal models still bear fundamental questions. The models emerge due to various applications like 
anomalous diffusion 
\cite{fractional_anomaDiff1,fractional_anomaDiff2}, 
peridynamics 
\cite{Silling_birthperidym}
%, peridym1, peridym2, peridym3, peridym4} 
or image processing 
\cite{imageproc1,imageproc2} and are diverse in their mathematical nature \cite{nonloc_vector_calc_large, bogdan2003censored}.

Their investigation is often accompanied with numerical experiments motivated by these applications. 
%2 description of the problem
In this work, we describe the discretization of nonlocal operators of the general form
$$-\mcL_\delta \ub (\xb) := ~ 2\int_{\base \cap B_\delta(\xb)} \Cb_\delta (\xb, \yb) \ub(\xb) -  \Cb_\delta (\yb, \xb) \ub(\yb) d\yb,$$
where the kernel $\Cb_\delta$ accounts for a finite range of nonlocal interactions determined by the horizon $\delta>0$.
The interaction neighborhood $B_\delta(\xb)$, i.e., the support of the kernel, is often given by a suitable approximation of the Euclidean or infinity norm ball.

The purpose of  \texttt{nlfem} is to compute numerical solutions to related boundary value problems which run at convenient speed for researchers. 
%3 what is the difficulty here
The discretization of these problems is achieved by a finite element approximation. The resulting variational framework comes at the price of a second integration layer compared to the operator in its strong form, which makes it more costly compared to classical differential operators. Further challenges in the implementation arise due to finite interaction horizons and singularities of the kernel.

%4 what has been done so far in this area
The development of efficient codes is an active field of research and there are various related implementations.
The foundation for operators related to the fractional Laplacian\footnote{Throughout this work we always refer to the fractional Laplacian in the integral form.} is given by
boundary element methods \cite{sauter_boundary_2011}. Based on these fundamentals, a Matlab implementation for a finite element approximation of the two--dimensional fractional Laplacian with infinite interaction is presented in \cite{acosta_FE_implementation}. Further advanced techniques to efficiently implement the fractional Laplacian and resulting dense stiffness matrices are developed in \cite{glusa_integralLap} and incorporated into the finite element code \textit{PyNucleus}\footnote{\href{https://github.com/sandialabs/PyNucleus}{https://github.com/sandialabs/PyNucleus}}.
The latter is a recommendable alternative to our code for
nonlocal diffusion problems with interaction neighborhoods based on the Euclidean norm ball. Specifically, it can handle large dense matrices which makes it well--suited for a large or even infinite horizon $\delta$.

Alternatives to kernel truncations have been investigated in \cite{aulisa2021efficient}.
Apart from this, various structure exploiting approaches, finite-difference schemes and meshfree methods have been used to compute numerical solutions. In particular meshfree methods currently dominate the practical application of nonlocal models for mechanics \cite{parks2011peridynamics, boys2021peripy}.
However, implementations of finite element methods for peridynamics have already been integrated into engineering software \cite{wu_dynamic_2018}.
For a general overview we refer the reader to the comprehensive review paper \cite{acta20} on numerical methods for nonlocal problems.\\

%5 what do we present (new)
    Our code \texttt{nlfem} assembles nonlocal operators on triangular meshes based on linear continuous (CG) or discontinuous Galerkin (DG) \textbf{ansatz spaces} in 2d and also offers this functionality in 1d and 3d for some generic settings. 

It allows the assembly of problems with given Neumann and Dirichlet \textbf{boundary data} \cite{SIREV, dipierroNeumann, foghem2022general}. An important detail is that the null space of systems corresponding to pure Neumann boundary data is exact up to machine precision, which allows an efficient evaluation of a pseudoinverse.

Concerning the \textbf{domain}, \texttt{nlfem} covers a variety of different configurations.
It handles nonlocal interactions in nonconvex or even disconnected domains where the intersection between the interaction neighborhood $B_\delta(\xb)$ and the domain can be disconnected. For example, this is of particular interest
in shape optimization with nonlocal operators \cite{Vollmann2019}, where the domain is modified iteratively. 

The \textbf{kernel} can be symmetric or nonsymmetric as well as scalar-- or matrix--valued. For the symmetric case our discretization of the weak form guarantees the symmetry of the stiffness matrix up to machine precision.
The code  can generically handle smooth kernels and it comes with quadrature rules for fractional type kernels as they are found in \cite{sauter_boundary_2011,acosta_FE_implementation}.
In addition to that, the kernel can depend on finite element labels which opens the door to a monolithic assembly of \textbf{interface problems} determined by spatially variable kernels. 
%
%\CV{Bezogen auf den letzten Satz: In welcher Weise unterschieden die sich? Kann man das vlt. in einem Satz beschreiben? Zudem: Bezieht sich das auf die Flexibilität, die wir für den Kern erlauben (aktueller Abschnitt)?}

%\gray{truncations}
 The \texttt{nlfem} code is most efficient for operators with \textbf{interaction horizons}, which are comparable to the mesh size, i.e., $h \leq \delta \leq C h$ for some $C \geq 1$. This relation is often used in the nonlocal mechanics setting; for example, in {\cite{Bobaru12,Parks08}}, the choices $\delta=3h$ and $\delta=4h$ are advocated, respectively.
In this regime the systems can be considered to be rather sparse.
A careful consideration of the quadrature and interpolation errors can allow smaller rations $\delta/h$ and increases the sparsity of the related systems. 

Our implementation is based on the extensive discussions on the errors incurred by various \textbf{approximations of the interaction neighborhood} and quadrature rules \cite{approxBall, Vollmann2019}, all of which are implemented here. We only discuss a small selection of interaction neighborhoods in this paper, such as two approximations of the Eucldiean norm ball, which provably do not deteriorate the finite element interpolation error \cite{approxBall}. Furthermore, the interaction neighborhood of a kernel is efficiently determined by a breadth--first traversal of finite elements throughout the assembly process, which avoids expensive preprocessing computations.

A fundamental advantage of finite element methods is that they can be considered to be \textbf{asymptotically compatible} in the sense of \cite{du_tian_locallim}.
Our code reproduces this property for the cases  $h \leq \delta$, $h \sim \delta \to 0$ if the implemented
interaction neighborhood does not induce geometric errors. For example, this is the case for the implemented infinity norm ball.

%\gray{solving}
For convenience the assembly is performed in multiple threads and the main routine, which is written in \textbf{C++}, comes with a user--friendly \textbf{Python} interface. When it comes to solving, we note that the stiffness matrix is returned in Compressed Sparse Row (CSR) format. Therefore, the user can apply any sparse solver accessible from Python and apply it to the stiffness matrix.\\

%6 organization of the paper
The remainder of this article is organized into two main Sections. First, in Section \ref{sec:fem} we give a precise formulation of the targeted problem class in \ref{subsec:problem_formulation} along with its finite element approximation in \ref{subsec:finiteDimensionalApproximation}. In its following subsections we then highlight discretization details that deserve special attention in a nonlocal framework. Second, in Section \ref{sec:numerics} we present various numerical examples, including diffusion and mechanics, and give a brief scaling study.

% PROBLEM FORMULATION
\section{Finite Element Approximation}\label{sec:fem}
\subsection{Problem Formulation}\label{subsec:problem_formulation}
Let $\base \subset \R^{\ddim}$ be a compact \textit{base domain}, and let $\Omega \subset \base$ be open in $\base$\footnote{More precisely, there exists an open set $\mathcal{O}\subset \R^d$, such that $\mathcal{O}\cap \base = \Omega$.}. We note that the case $\Omega = \base$ is allowable, and $\Omega$ is open in $\R^{\ddim}$ only if $\Omega \subset \interior(\base)$. We refer to $\Omega$ as \textit{domain}, although it is not necessarily open in $\R^{\ddim}$.
The complement of $\Omega$ in $\base$ is denoted by $\Omega_D := \base \setminus \Omega$. It can play the role of a nonlocal Dirichlet boundary in suitable settings. If $\base = \Omega$ the resulting nonlocal problem can be interpreted as Neumann type problem.
In Figure \ref{fig:dirichlet-domain} we present an exemplary configuration.
Let $\delta > 0$ be an \textit{interaction horizon}. We make the following assumption
about the kernel function.
\begin{assumption}
\label{assump:kernelFunction}
We assume that for the matrix--valued \emph{kernel function} $\Psi_\delta: \R^{\ddim} \times \R^{\ddim} \rightarrow \R^{{\ndim}\times {\ndim}}$ there exists 
an $s\in (0, 1)$ such that 
\begin{align}
    \Psi_\delta(\xb, \yb)\,
    \ind_{B_\delta(\xb)} (\yb)\,
    |\xb - \yb|^{{\ddim} + 2s} \label{ass:singularity}
\end{align}
is bounded.
\end{assumption}
This allows singularities at the origin and also includes smooth kernels such as the constant kernel. We assign to  the above kernel function the \textit{interaction neighborhood} $B_\delta(\xb)$, where
\begin{align}
     B_\delta(\xb) := \lbrace\yb \in \R^{\ddim}  ~|~ |\xb - \yb|_{\bullet} \leq \delta \rbrace,
\end{align}
for some norm $| \cdot  |_{\bullet}$ in $\R^d$. We denote the truncated \textit{kernel} by
\begin{equation}\label{eq:defiKernel}
    \Cb_\delta(\xb, \yb) := \Psi_\delta(\xb, \yb) \ind_{B_\delta(\xb)} (\yb).
\end{equation}
The linear \textit{nonlocal operator} under consideration acting on a function $\ub\colon \R^{\ddim}\to\R^{\ndim}$ is then given by
\begin{align}
    -\mcL_\delta \ub (\xb) := p.v. ~ 2\int_{\base} \Cb_\delta (\xb, \yb) \ub(\xb) -  \Cb_\delta (\yb, \xb) \ub(\yb) d\yb,
    \label{def:strongL}
\end{align}
where ``\textit{p.v.}'' denotes the Cauchy principal value\footnote{More precisely, $p.v. ~ \int_{\base} h(\xb,\yb) d\yb := \lim_{\epsilon\to 0^+} \int_{\base \setminus B_\epsilon(\xb)} h(\xb,\yb) d\yb$.}. By testing (\ref{def:strongL}) with $\vb \colon \base \to \R^{\ndim}$ where $\vb=0$ on $\Omega_D$,  we obtain the bilinear form
\begin{equation}\label{eq:biFormDef}
  \begin{aligned}
              A(\ub, \vb)&:=-\int_\Omega \vb(\xb)^\smalltop \mcL_\delta \ub (\xb) d\xb \\
    &=2 \int_{\Omega} \vb(\xb)^\smalltop  ~p.v. \int_{\base} \Cb_\delta(\xb, \yb)\ub(\xb)
    - \Cb_\delta(\yb, \xb) \ub(\yb) 
    d\yb  d \xb. 
  \end{aligned}  
\end{equation}
With $\vb=0$ on $\Omega_D$ and Fubini we observe that
\begin{align*}
&\int_{\Omega} \vb(\xb)^\smalltop  ~p.v. \int_{\base} \Cb_\delta(\xb, \yb)\ub(\xb)
    - \Cb_\delta(\yb, \xb) \ub(\yb) 
    d\yb  d \xb \\
=&\int_{\base} \vb(\xb)^\smalltop  ~p.v. \int_{\base} \Cb_\delta(\xb, \yb)\ub(\xb)
    - \Cb_\delta(\yb, \xb) \ub(\yb) 
    d\yb  d \xb \\
   =&-\int_{\base}  ~p.v. \int_{\base} \vb(\yb)^\smalltop \left(\Cb_\delta(\xb, \yb)\ub(\xb)
    - \Cb_\delta(\yb, \xb) \ub(\yb)\right)
    d\yb  d \xb.  
\end{align*}
So that the bilinear form can be written as
\begin{equation}
    \begin{aligned}
A(&\ub, \vb)\\ &= 
 \int_{\base} \int_{\base} (\vb(\xb)-\vb(\yb))^\smalltop  \left(\Cb_\delta(\xb, \yb)\ub(\xb)
    - \Cb_\delta(\yb, \xb) \ub(\yb) \right)
    d\yb  d \xb\\
    & =\int_{\base} \int_{\base  }\ind_{B_\delta(\xb)}(\yb) (\vb(\xb)-\vb(\yb))^\smalltop  \left(\Psi_\delta(\xb, \yb)\ub(\xb)
    - \Psi_\delta(\yb, \xb) \ub(\yb) \right)
    d\yb  d \xb.
\end{aligned}
\label{def:weakL} 
\end{equation}
Note that we exploited the symmetry of the exact indicator function $\ind_{B_\delta(\xb)} (\yb)$ to obtain this equality.
\subsection{Finite Dimensional Approximation}\label{subsec:finiteDimensionalApproximation}
Let $\mcT^h := \lbrace \mcE_\kdx \rbrace_{\kdx=1}^{\Kdx}$ denote a subdivision of $\base = \Omega \cup \Omega_D$ into polyhedral finite elements with nodes $\lbrace \xb_m \rbrace_{m=1}^\Mdx$. 
\begin{assumption}
We assume that $\Omega$ and $\Omega_D$ can be exactly covered by the subdivisions $\mcT^h_\Omega= \lbrace \mcE_\kdx \rbrace_{\kdx=1}^{K_\Omega}$ and $\mcT^h_D = \lbrace \mcE_k \rbrace_{\kdx=K_\Omega+1}^{\Kdx}$ with $\mcT^h = \mcT^h_\Omega \cup \mcT^h_D $, respectively, where $\mcT^h_D$ is possibly empty. Since we assume polyhedral elements, this implies that
\begin{align}
     \overline{\Omega} = \bigcup_{\kdx=1}^{K_\Omega} \overline{\mcE}_k
 \textrm{~~~and~~~}
     \overline{\Omega}_D = \bigcup_{\kdx=K_\Omega+1}^{\Kdx} \overline{\mcE}_k
\end{align}
are polyhedral domains.
\end{assumption}
In the case of $\Omega_D = \emptyset$ we have that $K_\Omega = \Kdx$ and $\mcT^h = \mcT^h_\Omega$. 
For convenience of notation we assume an ordering of the nodes such that $\lbrace \xb_m \rbrace_{m=1}^{M_\Omega} \subset \Omega$ and $\lbrace \xb_m \rbrace_{m=M_\Omega+1}^{\Mdx} \subset \overline{\Omega}_D$. This assumption is not made in the implementation of  \texttt{nlfem}.
\begin{figure}[ht]
    \begin{subfigure}{.5\textwidth}
	   \centering
 	   \includegraphics[scale=.15]{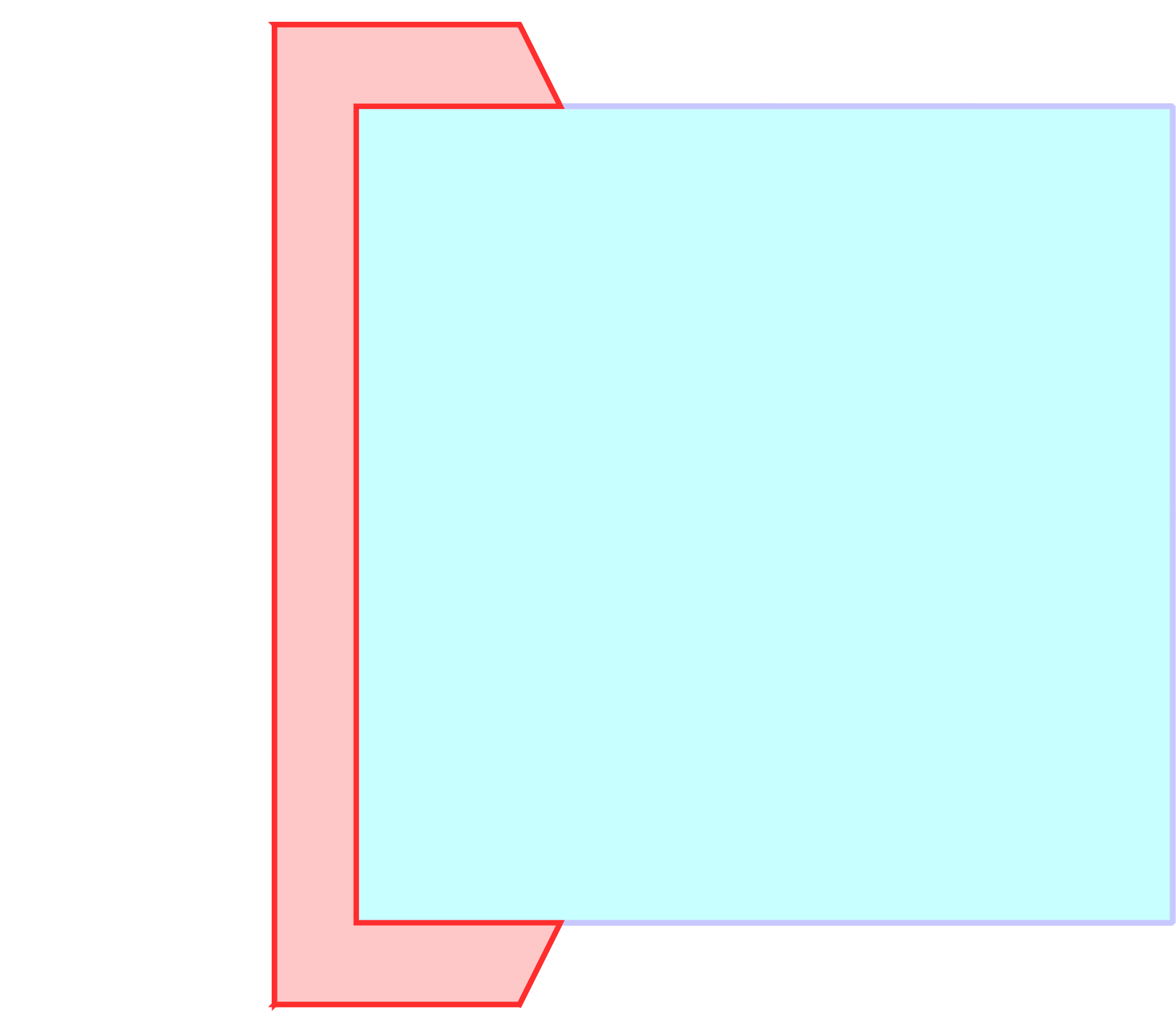}
 	   \caption{}
    \end{subfigure}
    \begin{subfigure}{.49\textwidth}
	    \centering
	    \includegraphics[scale=.15]{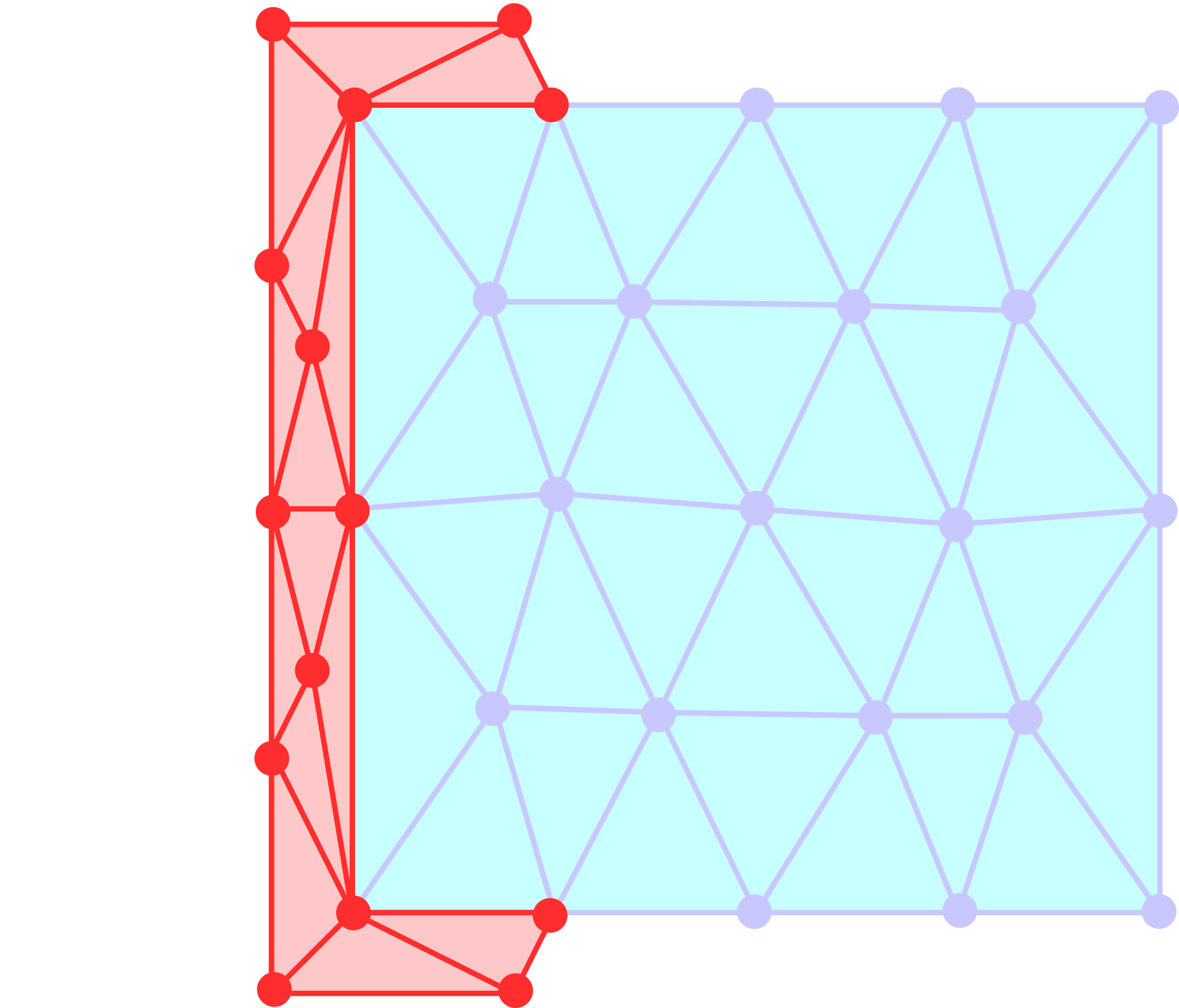}
	    \caption{}
    \end{subfigure}
    \caption{
	 (a) The compact base domain $\base$ contains the domains $\Omega$ (blue) and $\Omega_D$ (red), where the latter is possibly empty. The blue boundary belongs to $\Omega$ and the red boundary belongs to $\Omega_D$.
     (b) The red elements belong to the nonlocal Dirichlet boundary $\Omega_D$, while the blue elements belong to the domain $\Omega$. Note that the vertices on the red lines do not belong to the interior of $\Omega$.
    }
    \label{fig:dirichlet-domain}
\end{figure}
We implement piecewise--linear continuous and discontinuous ansatz functions $\lbrace \phi_{\jdx} \rbrace_{{\jdx}=1}^{\Jdx}$, where ${\Jdx} = nM$ in case of continuous and ${\Jdx} = {\ndim}({\ddim}+1)\Kdx$ in case of discontinuous ansatz functions.
Again for convenience we assume an ordering of the basis functions and define the corresponding finite dimensional subspaces
\begin{align}
    V^h(\base, \R^{\ndim}) := \spann ( \lbrace \phi_{\jdx} \rbrace_{{\jdx}=1}^{{\Jdx}} ), 
    \textrm{~~and~~}
    V^h_c(\base, \R^{\ndim}) := \spann ( \lbrace \phi_{\jdx} \rbrace_{{\jdx}=1}^{{\Jdx}_\Omega} ).
\end{align}
In case of continuous ansatz functions the unknown coefficients correspond to the nodes lying in the interior of $\Omega$ with respect to $\base$. See Figure \ref{fig:dirichlet-domain} for an illustration. Now, the evaluation of the bilinear form $A$ on $V^h(\base, \R^{\ndim}) \times V_c^h(\base, \R^{\ndim})$ can be written as sum over the finite elements, i.e.,
\begin{equation}
\begin{aligned}
    & A(\phi_{\jdx}, \phi_{\idx})= \\
    & \sum_{\kdx=1}^{\Kdx} \sum_{\ldx=1}^{\Kdx} 
     \left[   \int_{\mcE_k}   \int_{\mcE_{\ldx}} (\phi_{\idx}(\xb)-\phi_{\idx}(\yb))^\smalltop \left(
     \Cb_\delta(\xb, \yb)\phi_{\jdx}(\xb) 
     - \Cb_\delta(\yb, \xb) \phi_{\jdx}(\yb)\right)
     d\yb d \xb \right].
\end{aligned}
\label{eq:biFormOverElements}%
\end{equation}
Since the kernel  $\Cb_\delta$ may exhibit a truncation on some pairs $(\mcE_k, \mcE_{\ldx})$ we need an appropriate approximation of its support.
Apart from a modeling aspect, a numerical advantage of the kernel truncation is that it produces sparse systems on sufficiently coarse grids. However, this can only be exploited if the numerical evaluation of the interaction neighborhood $B_\delta(\xb)$ does not deteriorate the interpolation error of the finite element subspace. While the approximation of the infinity or Manhattan norm balls does not introduce a geometric error, the one for curved neighborhoods, like the Euclidean norm ball, does.

\begin{figure}[ht]
    \begin{subfigure}{.5\textwidth}
    	    \centering
    		\includegraphics[scale=.25]{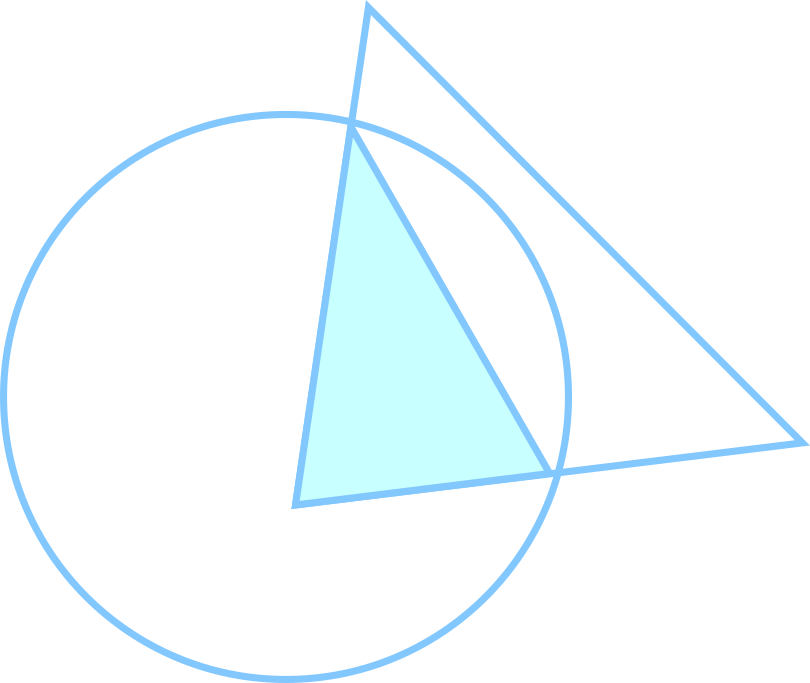}
    		\caption{}
   	\end{subfigure}
	\begin{subfigure}{.49\textwidth}
		\centering
        \includegraphics[scale=.25]{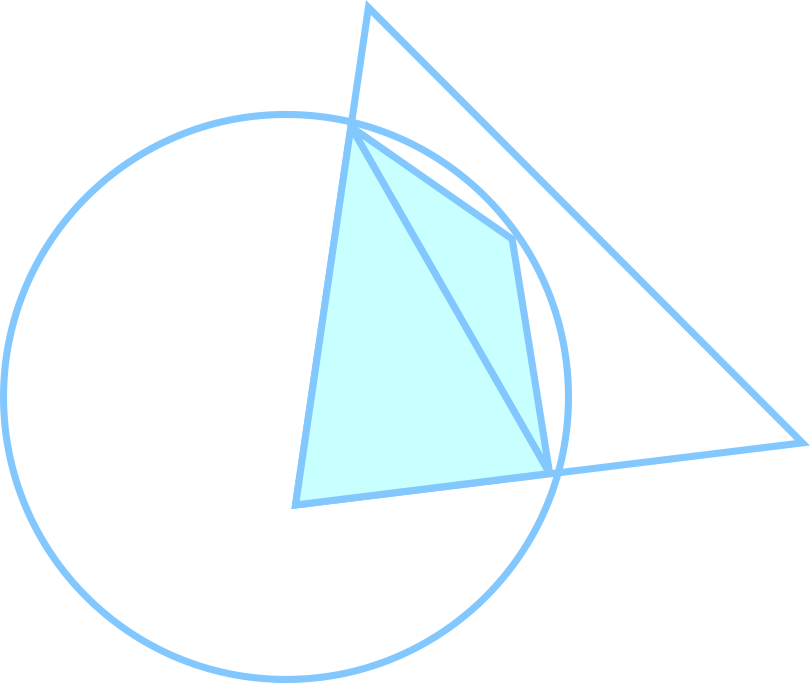}
        \caption{}
	\end{subfigure}

    \caption{
        If a finite element triangle is only partially covered by an Euclidean ball
        the intersection contains circular caps. The \NoCaps ball (a) omits this cap, whereas the
        \ApproximateCaps ball (b) retriangulates the whole intersection.
    }
    \label{fig:intersection-caps}
\end{figure}

In the following, we describe two major examples of the implemented ball approximations for the Euclidean ball 
$$B_{\delta}^2(\xb) := \lbrace \yb\in \R^d ~ | ~ |\xb - \yb|_2 \leq \delta \rbrace,$$ 
which we call the \ApproximateCaps and the \NoCaps approximations. Both are based on the given finite element mesh $\mcT^h$. By ``cap'' we mean the circular segments that arise when a finite element triangle is only partially covered by the Euclidean ball; see Figure \ref{fig:intersection-caps}. Among others, these ball approximations are investigated in \cite{approxBall}. 
\begin{definition}[\NoCaps Ball]
\label{def:noCapsBall}
For $\xb \in \base$ the \NoCaps ball approximation is defined as the convex hull of the intersection of the boundary $\partial B^2_\delta(\xb)$ of the Euclidean ball and the boundaries of the elements, i.e.,
$$
B_{\delta}^{ncp}(\xb) := \conv\left( 
    \bigcup_{\mcE_\ldx \in \mcT^h} \partial\mcE_\ldx\cap \partial B_{\delta}^2(\xb)
\right),
$$
where the dependency on the mesh size $h$ is omitted in the notation.
\end{definition}
 For ${\ddim}=2$, one can show that the area of the symmetric difference of $B^2_{\delta}(\xb)$ and $B_{\delta}^{ncp}(\xb)$
 is provably of order $\mcO(h^2)$ on a given mesh of size $h$. This is error commensurate with respect to the interpolation error of a linear finite element ansatz space; details for the latter two statements can be found in \cite{approxBall}.
The convex hull omits circular caps which appear in the intersection of some elements $\mcE_{\ldx} \cap B_{\delta}^2(\xb)$ and have a significant size on coarse grids. Therefore, by adding additional points on the center of possible caps the geometric error can be reduced even further while still being in $\mcO(h^2)$ for $\ddim=2$; see again Figure \ref{fig:intersection-caps}. The maximum number of caps for a single intersection $\mcE_{\ldx} \cap B_{\delta}^2(\xb)$ is three.

While the results in \cite{approxBall} are derived for a fixed horizon $\delta$, related investigations for the local limit ($\delta \to 0$) with polygonal ball approximations can be found in \cite{Du2021:locallimitApproxBall}.

\begin{definition}[\ApproximateCaps Ball]
\label{def:approxCapsBall}
Let $\xb\in  \base$. We denote the points on the cap center of each nonempty intersection $\mcE_{\ldx} \cap \partial B_{\delta}^2(\xb)$ by $\yb_{\ldx}$. 
The \ApproximateCaps ball is then defined by
\begin{align*}
    B_{\delta}^{acp}&(\xb) := \\
    &\conv
    \left(
        B_{\delta}^{ncp}(\xb)
    \right. 
     \cup 
      \left.
         \,\lbrace 
            \yb_{\ldx} \, | \, \yb_{\ldx} \textrm{ cap center of } \mcE_{\ldx} \cap \partial B_{\delta}^2(\xb) \textrm{ for some } \mcE_{\ldx} \in \mcT^h
        \rbrace 
     \right).
\end{align*}    
\end{definition}

Exact quadrature rules for circular caps can be found in \cite{fies_algebraic_nodate}. Here however, the quadrature points have to be computed during run time, as the rules depend on the geometry of the cap for higher quadrature orders. Therefore, we do not consider these exact rules in \texttt{nlfem}.

In addition to these approximate balls, we next also introduce the \InfinityNorm ball.
\begin{definition}[Infinity Norm Ball]
    \label{def:inftyBall}
    For $\xb \in  \base$ the \InfinityNorm ball is defined by 
    $$
        B_{\delta}^{\infty}(\xb) := \lbrace \yb \in\R^d~ | ~ |\xb - \yb|_\infty \leq \delta \rbrace.
    $$
\end{definition}

Proofs for the convergence of the nonlocal Dirichlet problem to the classical Dirichlet problem with corresponding scaling constants for various kernel functions can be found in \cite{Vollmann2019}. Since the \InfinityNorm ball is implemented exactly, it allows numerical tests of the expected asymptotic compatibility of the discretization scheme \cite{du_tian_locallim}.

We finally note that the indicator function based on any implemented truncation may lack in symmetry. More precisely, there might exist $\xb$ and $\yb$, for which
\begin{align*}
    \ind_{B_{\delta}^\#(\xb)}(\yb) \neq \ind_{B_{\delta}^\#(\yb)}(\xb),
\end{align*}
where $B_{\delta}^\#(\xb)$, $\#\in\{ncp,acp,\infty\}$, represents one of the implemented truncations.
This artefact stems from the ball approximation itself in case of the \NoCaps and \ApproximateCaps ball, but can also be caused by the quadrature; see Remark \ref{rem:assymmetryDueToQuadrature} below.

We define the integrand
\begin{equation}\label{def:integrand}
    \Phi_{\idx {\jdx} }(\xb, \yb) :=  ( \phi_{\idx}(\xb) - \phi_{\idx}(\yb))^\smalltop   
    (
        \Psi_\delta(\xb, \yb) \phi_{\jdx}(\xb)
        - \Psi_\delta(\yb, \xb) \phi_{\jdx}(\yb) 
    )
\end{equation}
and, based on the ball approximation and \eqref{def:weakL}, the \textit{approximate bilinear form}
\begin{equation}\label{eq:biFormApprox}%
\begin{aligned}
     A_h^{\#}(\phi_{\jdx}, \phi_{\idx})
     &:= \int_{\base}\int_{\base}\ind_{B_{\delta}^\#(\xb)}(\yb)
	\Phi_{ \idx {\jdx}}(\xb, \yb)d\yb  d \xb\\
	&=
     \sum_{\kdx=1}^{\Kdx} \sum_{\ldx=1}^{\Kdx} 
     \int_{\mcE_k}   \int_{\mcE_{\ldx}} 
     \ind_{B_{\delta}^\#}(\xb, \yb) \Phi_{\idx {\jdx} }(\xb, \yb) d\yb d \xb .
\end{aligned}
\end{equation}
By Fubini's integration theorem, for sufficiently smooth ansatz functions $\phi_{\jdx} \in V^h(\base, \R^{\ndim})$ and $\phi_{\idx}  \in V_c^h(\base, \R^{\ndim})$, the approximate bilinear form defined in \eqref{eq:biFormApprox} can be written as
\begin{align*}
&A_h^{\#}(\phi_{\jdx}, \phi_{\idx}) \\
    &= \int_{\base}\int_{\base}\ind_{B_{\delta}^\#(\xb)}(\yb)
	\Phi_{ \idx {\jdx}}(\xb, \yb)d\yb  d \xb\\
    &=\int_{\Omega} ~p.v. \int_{\base}\ind_{B_{\delta}^\#(\xb)}(\yb) \phi_{\idx}(\xb)^\smalltop 
	\left(
	\Psi_\delta(\xb, \yb) \phi_{\jdx}(\xb)  
	- \Psi_\delta(\yb, \xb) \phi_{\jdx}(\yb) 
	\right) d\yb ~ d \xb  \\
	 &~~~-\int_{\base} ~p.v. \int_{\Omega}\ind_{B_{\delta}^\#(\xb)}(\yb) \phi_{\idx}(\yb)^\smalltop 
	\left(
	\Psi_\delta(\xb, \yb) \phi_{\jdx}(\xb)  
	- \Psi_\delta(\yb, \xb) \phi_{\jdx}(\yb) 
	\right) d\yb ~ d \xb  \\
    &=
    2 \int_{\Omega} \phi_{\idx}(\xb)^\smalltop   ~p.v. \int_{\base}   \ind^S_{B_{\delta}^\#} (\xb, \yb) 
    \left(
       \Psi_\delta(\xb, \yb) \phi_{\jdx}(\xb)  
        - \Psi_\delta(\yb, \xb) \phi_{\jdx}(\yb) 
   \right) d\yb ~ d \xb ,
\end{align*}
where
\begin{equation} \label{eq:symmetrifiedIndicator}
    \ind^S_{B_{\delta}^\#}(\xb, \yb):=\frac{1}{2}\left( \ind_{B_{\delta}^\#(\xb)}(\yb) + \ind_{B_{\delta}^\#(\yb)}(\xb) \right).
\end{equation}
In view of \eqref{eq:biFormDef}, this shows that the approximate bilinear form $A_h^{\#}$ can be interpreted as the discretization of the operator $-\mcL_\delta$ based on the symmetrified approximate indicator function $\ind^S_{B_{\delta}^\#}(\xb, \yb)$ instead of $\ind_{B_{\delta(\xb)}}(\yb)$ in the strong form. However, defining the approximate bilinear form as in \eqref{eq:biFormApprox} guarantees the symmetry of the stiffness matrix $(A_h^{\#}(\phi_{\jdx}, \phi_{\idx}) )_{\idx,\jdx}$; see also Remark \ref{rem:assymmetryDueToQuadrature} below.

In the following we discuss the evaluation of the local contributions
\begin{align}
    \localStiffness{\kdx}{\ldx}{\jdx}{\idx} := \int_{\mcE_k}\int_{\mcE_{\ldx}}\ind_{B_{\delta}^\#(\xb)}(\yb)\Phi_{\idx {\jdx}}(\xb, \yb)d\yb ~ d \xb
    \label{def:localStiffness}
\end{align}
to the $(i,j)$-th entry of the stiffness matrix.

\subsection{Population of the Stiffness Matrix}\label{subsec:populationStiffnessMatrix}
The assembly algorithm iterates over all pairs of elements and then adds the local contributions \eqref{def:localStiffness} to the stiffness matrix. We call a pair of elements $(\mcE_{\kdx}, \mcE_{\ldx})$ \emph{intersecting} if
$$
    \overline{\mcE_k} \cap \overline{\mcE_{\ldx}} \neq \emptyset.
$$
If a pair is not intersecting, we call it \emph{disjoint}. 
For the evaluation of the contributions it is only important how the kernel behaves on a fixed pair of elements. If for example the kernel exhibits a singularity at the origin and the pair is disjoint, the singularity does not occur. Also, if a pair of elements $(\mcE_{\kdx}, \mcE_{\ldx})$ fulfills that
\begin{align*}
    \mcE_{\ldx} \subset B_\delta(\xb) \textrm{ for all } \xb \in \mcE_k,
\end{align*}
the truncation does not come into play. We therefore distinguish two cases. In the first case the kernel has a singularity with $s \leq 0.5$ or the pair of elements is disjoint. In the second case the kernel has a singularity with $s>0.5$ and the elements are intersecting.

\subsubsection{Disjoint Pairs or Kernels with $s\leq 0.5$}
If $s \leq 0.5$ in \eqref{ass:singularity} or if the pair $(\mcE_{\kdx}, \mcE_{\ldx})$ is disjoint, the local contributions (\ref{def:localStiffness}) can be factored out and computed separately, i.e.,
\begin{equation}\label{eq:case-one-factoredOut}
    \begin{aligned}
    \localStiffness{\kdx}{\ldx}{\jdx}{\idx}= &
    % LOCAL
    \int_{\mcE_{\kdx}} \int_{\mcE_{\ldx}}\ind_{B_{\delta}^\#(\xb)}(\yb) \phi_{\idx}(\xb)^\smalltop \Psi_\delta(\xb, \yb) \phi_{\jdx}(\xb)   d\yb ~ d \xb  \\
    % NON LOCAL    
    &-\int_{\mcE_{\kdx}} \int_{\mcE_{\ldx}}\ind_{B_{\delta}^\#(\xb)}(\yb)
        \phi_{\idx}(\xb)^\smalltop \Psi_\delta(\yb, \xb) \phi_{\jdx}(\yb)   d\yb ~ d \xb  \\
    % LOCAL PRIME
    &+\int_{\mcE_{\kdx}} \int_{\mcE_{\ldx}}  \ind_{B_{\delta}^\#(\xb)}(\yb)
    \phi_{\idx}(\yb)^\smalltop \Psi_\delta(\yb, \xb) \phi_{\jdx}(\yb)  d\yb ~ d \xb  \\
    % NONLCAL PRIME
    &-\int_{\mcE_{\kdx}} \int_{\mcE_{\ldx}} \ind_{B_{\delta}^\#(\xb)}(\yb) \phi_{\idx}(\yb)^\smalltop \Psi_\delta(\xb, \yb) \phi_{\jdx}(\xb)  d\yb ~ d \xb. \\
    =: & 
    \termLocal{\kdx}{\ldx}{\jdx}{\idx} 
    + \termNonloc{\kdx}{\ldx}{\jdx}{\idx} \\
    &+ \termLocalPrime{\kdx}{\ldx}{\jdx}{\idx} 
    + \termNonlocPrime{\kdx}{\ldx}{\jdx}{\idx}.
\end{aligned}
\end{equation}
Therefore, \texttt{nlfem} allows continuous and discontinuous ansatz functions if $s\leq 0.5$. Note, that the fractional Sobolev spaces $H^s(\base)$ contain functions with jump discontinuities for $0\leq s\leq 0.5$. 

The expression  $\termLocal{\kdx}{\ldx}{\jdx}{\idx}$ is nonzero only if the element $\mcE_\kdx$ lies in the support of $\phi_\idx$ and $\phi_\jdx$. Similarly, the contribution of $\termLocalPrime{\kdx}{\ldx}{\jdx}{\idx}$ is linked to the element $\mcE_\ldx$. The term $\termNonloc{\kdx}{\ldx}{\jdx}{\idx}$ is nonzero only if $\phi_\idx$ has its support on $\mcE_\kdx$ and $\phi_\jdx$ on $\mcE_\ldx$, where the converse holds for $\termNonlocPrime{\kdx}{\ldx}{\jdx}{\idx}$.  
That way we derive the indices of the basis functions corresponding to the pair $(\mcE_\kdx, \mcE_\ldx)$ in the stiffness matrix.  We note that the same indices of basis functions occur again for the pair $(\mcE_\ldx, \mcE_\kdx)$, but the contributions to the stiffness matrix are not identical as the truncation $\ind_{B_{\delta}^\#(\xb)}(\yb)$ is not symmetric. Ultimately, the contribution of the two pairs $(\mcE_\kdx,\mcE_\ldx)$ and $(\mcE_\ldx,\mcE_\kdx)$  together lead to the symmetrified truncation \eqref{eq:symmetrifiedIndicator}.

\subsubsection{Kernels with $s > 0.5$ on Intersecting Pairs}
If a pair of elements is intersecting and $s > 0.5$ the separation of integrands in \eqref{eq:case-one-factoredOut} is not admissible. Therefore,  \texttt{nlfem} is restricted to continuous ansatz functions if $s > 0.5$. As the pair is intersecting, the elements are either vertex touching, edge touching or identical. 
Therefore the number of ansatz functions to be considered in these cases is 5, 4 or 3, respectively. 
If we denote the vertices by $\jdx_1, \dots, \jdx_\ldx$ for $\ldx = 5,4,3$ we obtain 25, 16, or 9 pairs of ansatz functions $(\phi_{\jdx_\nu}, \phi_{\jdx_{\nu'}})$ for $\nu, \nu' = 1, \dots, \ldx$ which yield nonzero contributions to the local stiffness matrix $(\localStiffness{\kdx}{\ldx}{\jdx_\nu}{\jdx_{\nu'}})_{1\leq \nu,\nu'\leq \ell}$.

\subsection{Quadrature}\label{subsec:quadrature}
The quadrature rules need to work for kernels with truncations and singularites. However, the implemented quadrature rules in \texttt{nlfem} do not take account of both at the same time on a fixed pair of elements. We therefore require the following assumption.
\begin{assumption}
The quadrature rules for singular kernels assume that for all intersecting pairs of elements $(\mcE_k, \mcE_\ell)$ it holds that $\mcE_{\ldx} \subset B_\delta(\xb)$ for all $\xb \in \mcE_k$.
\label{assump:fullyInteract}
\end{assumption}
\begin{rem}
    Specific quadrature rules for singular kernels are required for sufficiently strong singularities only. Kernels like the peridynamic kernel (\ref{ker:peridynamics}) can be integrated by simply 
    technically avoiding zero-divisions and we ignore those cases in our discussion of singular kernels.
\end{rem}
Given Assumption \ref{assump:fullyInteract} and the fact that disjoint pairs do not require a treatment with regularizing integral transforms we can evaluate all contributions as given in Algorithm \ref{alg:populationStiffnessMatrix}.
\begin{algorithm}
	\caption{Evaluate $\localStiffnessMatrix{\kdx}{\ldx}$ in \eqref{def:localStiffness}} 
	\label{alg:populationStiffnessMatrix}
	\begin{algorithmic}[1]
	    \If {$\mcE_{\kdx}, \mcE_{\ldx}$ intersect \textbf{and} the kernel is singular}
	        \State Evaluate $\localStiffnessMatrix{\kdx}{\ldx}$ via \textbf{quadrature for singular kernels}; 
	        \State see Section \ref{subsec:quad-singular} below
	    \Else 
	        \State Evaluate $\localStiffnessMatrix{\kdx}{\ldx}$ via \textbf{quadrature for kernel truncations};
	        \State see Section \ref{subsec:quad-truncation} below
	    \EndIf
	\end{algorithmic} 
\end{algorithm}

In any case, the quadrature is performed by pulling back the domain of integration to a reference domain 
\begin{align*}
    \widehat{\mcE} \times \widehat{\mcE},~~ \textrm{ where }~~ \widehat{\mcE}  := \lbrace \xbh\in \R^d ~|~ \xbh \geq 0, \sum_{\iota=1}^{{\ddim}} \xbh_\iota \leq 1\rbrace.
\end{align*}
The affine linear mapping $\chi_k: \widehat{\mcE} \rightarrow \mcE$ from the reference to a physical element allows to define the pullback
\begin{align*}
    \widehat{\Phi}_{\kdx\ldx,  \idx {\jdx}}(\xbh, \ybh) := 
    \Phi_{ \idx {\jdx}}(\toPhys{\kdx}{\xbh}, \toPhys{\ldx}{\ybh}),
\end{align*}
the inner integral
\begin{align}
        \widehat{\mcK}_{\kdx\ldx, \idx {\jdx}}^\#(\xbh) := \int_{\widehat{\mcE}}
    \ind_{B_{\delta}^\#} (\toPhys{\kdx}{\xbh}, \toPhys{\ldx}{\ybh}) \widehat{\Phi}_{\kdx\ldx,  \idx {\jdx}}(\xbh,\ybh) d\ybh
    \label{def:innerReferenceIntegrand}
\end{align}
and the local contribution to the $({\idx},\jdx)$-th entry of the stiffness matrix 
\begin{align}
        \widehat{A}_{\kdx\ldx, \idx {\jdx}}^\# :=
        \int_{\widehat{\mcE}}
        \widehat{\mcK}_{\kdx\ldx, \idx {\jdx}}^\#(\xbh)
    d\xbh
    \label{def:referenceLocalStiffness}
\end{align}
on the reference element.
\subsubsection{Quadrature for Kernel Truncations} \label{subsec:quad-truncation}
For some pairs $(\mcE_k, \mcE_{\ldx})$ we find that $\mcE_{\ldx}$ is only partially covered by the interaction neighborhood $B_{\delta}^\# (\xb)$ for some $\xb \in \mcE_k$, so that $\mcE_{\ldx} \cap B_{\delta}^\# (\xb) \subsetneq \mcE_{\ldx}$. %for those $\xb \in \mcE_k$. 
In this case, the ball approximations $B_{\delta}^{ncp}(\xb)$, $B_{\delta}^{acp}(\xb)$ as well as the ball $B_{\delta}^{\infty}(\xb)$ require a retriangulation of the integration domain $\mcE_{\ldx} \cap B_{\delta}^\# (\xb)$. 
We denote the set of elements which result from such a retriangulation by
\begin{align*}
    \mcT_{h, \ldx}^\#(\xb) := \lbrace \tilde{\mcE}_\ldx \rbrace_{\ldx=1}^{L_{\ldx}},
\end{align*}
so that $\bigcup \mcT_{h, \ldx}^\#(\xb) = \mcE_{\ldx} \cap B_{\delta}^\# (\xb)$. Let $\lbrace \ybh_{\qdx}, \quadyq \rbrace_{\qdx=1}^{\Qdx}$ denote a quadrature rule on the reference element $\widehat{\mcE}$, and let $\xbh \in \widehat{\mcE}$ be some reference point.
Then the fully-discrete inner integral from (\ref{def:innerReferenceIntegrand}) reads as
\begin{equation}\label{eq:discrete-inner}
 \widehat{\mcK}_{\kdx\ldx, \idx {\jdx}}^\#(\xbh)
 \approx
    \sum_{\tilde{\ldx}= 1}^{L_{\ldx}} \sum_{ \qdx = 1}^{\Qdx}
     \widehat{\Phi}_{\kdx \tilde{\ldx},  \idx {\jdx}}(\xbh,\ybh_{\qdx}) \quadyq.
\end{equation}
In finite element implementations the function values of ansatz functions are usually pre-computed at the quadrature points and stored. However, if a retriangulation is necessary the physical coordinates of the quadrature points 
$\toPhys{\tilde{\ldx}}{\ybh_{\qdx}}$ for $\qdx = 1,\dots,\Qdx$ on some element 
$\tilde{\mcE}_{\tilde{\ldx}} \in \mcT_{h, \ldx}^\#(\xb)$ are known at runtime only and the ansatz functions are evaluated at the corresponding points.

Now let $\lbrace \xbh_{ \pdx}, \quadxp \rbrace_{\pdx=1}^{\Pdx}$ be a quadrature rule for the reference element of the outer integral. Then with \eqref{eq:discrete-inner} the discretized version of the local contribution to the 
$(\idx, \jdx)$--th entry of the stiffness matrix \eqref{def:referenceLocalStiffness} is obtained by 
\begin{align*}
\widehat{A}_{\kdx\ldx, \idx {\jdx}}^\#
\approx
    \sum_{\pdx=1}^\Pdx \left(
    \sum_{\tilde{\ldx} = 1}^{L_{\ldx}} \sum_{ \qdx = 1}^{\Qdx}
    \widehat{\Phi}_{\kdx\tilde{\ldx},  \idx {\jdx}}(\xbh_{\pdx},\ybh_{\qdx}) \quadyq\right) \quadxp .
\end{align*}
\subsubsection{Quadrature for Singularities} \label{subsec:quad-singular}
Some kernel functions exhibit a singularity at the origin, which lies in the integration domain of 
$\widehat{\mcK}^\#_{\kdx\ldx, \idx {\jdx} }(\widehat{\xb})$
whenever
the pair $(\mcE_{\kdx}, \mcE_{\ldx})$ intersects.  
We therefore require regularizing integral transforms \cite{sauter_boundary_2011}. Assumption \ref{assump:fullyInteract} allows to ignore possible truncations in \eqref{def:referenceLocalStiffness} and to simply evaluate
\begin{align}
      \widehat{A}_{\kdx\ldx, \idx {\jdx}}^\#
      %\widehat{A}_{\kdx\ldx, \idx {\jdx} }^\# 
      =  \int_{\widehat{\mcE}} \int_{\widehat{\mcE}}
        \widehat{\Phi}_{\kdx\ldx, \idx  {\jdx}}(\xbh, \ybh) d\ybh d\xbh.
    \label{def:innerCloseReferenceIntegral}
\end{align}
If Assumption \ref{assump:fullyInteract} is violated, then the support of the kernel is overestimated by \eqref{def:innerCloseReferenceIntegral}.
Intersecting pairs $(\mcE_{\ldx},\mcE_k)$ can be vertex touching, edge touching or identical. For each of those cases we apply integral transforms, which again 
pull back subsets of the integration domain $\widehat{\mcE} \times \widehat{\mcE}$ to the unit cube $(0,1)^4$. 
The transformations are well established and applied for example in the field of boundary element methods. Details can be found in \cite{sauter_boundary_2011, acosta_FE_implementation}.
Note again, that in the case of singular kernels with $s>0.5$ the implemented routines cannot evaluate the expression in (\ref{def:innerCloseReferenceIntegral}) for discontinuous basis functions.

\begin{figure}[ht]
    \centering
    \includegraphics[scale=.30]{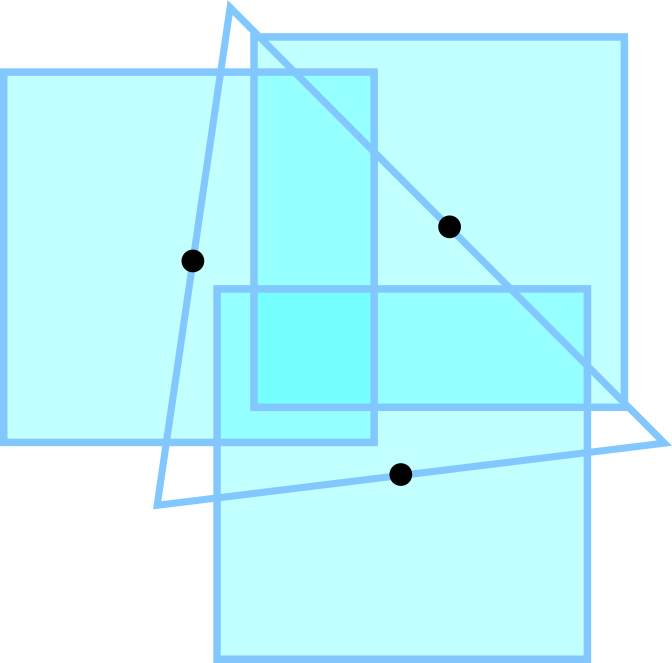}
    \caption{The Figure shows the approximation of the interaction domain of a single triangular element by
    interaction neighborhoods of three quadrature points for the \InfinityNorm ball $B_\delta^\infty$.}
    \label{fig:retriangulate-infty}
\end{figure}

\begin{rem}[Asymmetry due to quadrature] 
\label{rem:assymmetryDueToQuadrature}
    We have mentioned in Section \ref{subsec:finiteDimensionalApproximation} that the approximate indicator function based on the \ApproximateCaps and \NoCaps balls may generally lack in symmetry, and therefore would lead to nonsymmetric stiffness matrices if one used representation \eqref{eq:biFormDef} of the bilinear form. Also truncations invoked by the Infinity ball $B_\delta^\infty(\xb)$, which can be implemented exactly, can be nonsymmetric on sufficiently irregular grids. These two observations holds true \textit{independent} of the symmetry of the kernel function $\Psi_\delta(\xb,\yb)$. However, it is a desirable feature that the symmetry of the kernel, i.e., the self-adjointness of the operator, is transported through the discretization process. With other words, for a symmetric kernel we expect a symmetric stiffness matrix.
     This behaviour is more intricate and related to the approximation of the interaction domain of a finite element by the union of interaction neighborhoods of quadrature points, see Figure \ref{fig:retriangulate-infty}. Thus, we use by default representation \eqref{eq:biFormApprox} instead, which as stated above corresponds in strong form to the operator $-\mcL_\delta$ based on the approximate indicator function \eqref{eq:symmetrifiedIndicator}. Thereby, we guarantee the stiffness matrix to be symmetric up to machine precision for any symmetric kernel function. 
     Also, an important consequence of the symmetrification is that
     the null space of the stiffness matrix related to pure Neumann type problems contains the constant vectors up to machine precision. 
     This allows the application of projected Krylov subspace methods or
     rank-1 corrections to efficiently evaluate a pseudoinverse of the stiffness matrix.
\end{rem}

We summarize the different cases for the quadrature of kernels in Table \ref{tab:quadratureCases}. The table distinguishes between different degrees of singularity $s$ in case of $\ddim=2$ as well as intersecting and disjoint pairs of elements. 
When $s > 0$ we apply regularizing integral transforms represented by $\duffyQuadRule$ in the table. This quadrature rule is applied for intersecting pairs and requires Assumption \ref{assump:fullyInteract}. 
Whenever $s \leq 0.5$ we allow for discontinuous and continuous ansatz functions. When $s > 0.5$ only continuous ansatz functions are provided. For more details we refer to \cite{glusa_integralLap}.
If $-1 < s \leq 0$ the singularity might be so weak, that it suffices to simply avoid zero-divisions on intersecting pairs, which is represented by $\noZeroDivisionQuadRule$ in the table. However, the quadrature rule $\duffyQuadRule$ can also be applied.
If $s\leq -1$ there is no singularity at all and symbol $\commonQuadRule$ represents a standard quadrature rule.
The rules, $\commonQuadRule$ and $\noZeroDivisionQuadRule$, respect the kernel truncation, i.e. they do not require Assumption \ref{assump:fullyInteract}. 

\begin{table}[ht]
    \centering
    \begin{tabular}{|l|c|c|c|c|c|}
    \hline
                        & $s \leq -1$   & $-1 < s \leq 0$   & $0 < s \leq 0.5$  & $0.5 < s < 1$ \\ 
    \hline
         Intersecting   & $\commonQuadRule $             & $\noZeroDivisionQuadRule$ or $\duffyQuadRule$       & $\duffyQuadRule$      & $\duffyQuadRule $   \\
         Disjoint       & $\commonQuadRule$   & $\commonQuadRule$   & $\commonQuadRule$        &  $\commonQuadRule$                 \\ 
    \hline
                        & \multicolumn{3}{|c|}{CG and DG}                              &         CG      \\
    \hline
    \end{tabular}
{\centering\small~\\[0.2cm]
\begin{tabular}{ll}
$\commonQuadRule$&Standard quadrature rule on $\widehat{\mcE} \times \widehat{\mcE}$\\
$\noZeroDivisionQuadRule$& Quadrature rule on $\widehat{\mcE} \times \widehat{\mcE}$ that avoids $\xbh_{\pdx}=\ybh_{\qdx}$\\
$\duffyQuadRule$& Quadrature rule on $(0,1)^4$ after regularizing integral transforms
\end{tabular}
}
    \caption{The table gives an overview of the different quadrature rules which can be applied for dimension $\ddim=2$ in  \texttt{nlfem}.}
    \label{tab:quadratureCases}
\end{table}

\subsection{Traversal of the interaction neighborhood}\label{subsec:traversal}
The local contributions $\widehat{A}_{\kdx\ldx, \idx {\jdx}}^\#$ to the stiffness matrix, defined in \eqref{def:localStiffness}, can be nonzero even for pairs of remote finite elements $(\mcE_{\kdx}, \mcE_{\ldx})$, which we then refer to as \textit{interacting} elements in the following.
Unstructured meshes and the various supports of kernel functions
call for a flexible routine to identify interacting elements.
The identification can be accomplished if the assembly follows a breadth--first search. 
In order to describe the algorithm we define the \textit{adjacency graph} $\mbT_{adj} := \left(\mcT^h ,E_{adj} \right)$ of the finite element mesh with vertices $\mcT^h$ and edges
$$E_{adj} := \lbrace (\mcE_k, \mcE_{\ldx}) \in \mcT^h\times \mcT^h \, | 
\, \overline{\mcE_k} \cap \overline{\mcE_{\ldx}} \neq \emptyset\rbrace.$$
The graph $\mbT_{adj}$ can be understood as the dual graph of the finite element mesh, and its vertices are therefore given by the elements.
We additionally define the \textit{interaction graph} $\mbT_S := \left(\mcT^h, E_S\right)$ with vertices $\mcT^h$ and edges
\begin{align}
\label{graph:interactionEdges}
    E_S := \left\lbrace \left. (\mcE_k, \mcE_{\ldx}) \in \mcT^h\times \mcT^h \, \right| \, 
        \widehat{A}_{\kdx\ldx, \idx {\jdx}}^\# \neq 0 
        ~\text{for some}~
        \phi_\jdx,\phi_\idx \in V^h(\base, \R^n)
    \right\rbrace.
\end{align}
It is clear that $\mbT_{adj}$ is a spanning subgraph of $\mbT_S$ for any positive interaction horizon $\delta$. The adjacency graph $\mbT_{adj}$ can be computed and stored efficiently, while the interaction graph $\mbT_S$ exhibits storage requirements which are comparable to the full stiffness matrix.
We also note that the breadth--first search described below allows to naturally identify intersecting pairs $(\mcE_k, \mcE_\ell)$, as they only occur in the first layer of the traversal. This important built--in feature is used to identify the intersection cases for element pairs mentioned in Section \ref{subsec:populationStiffnessMatrix}.

Furthermore, for a fixed element $\mcE_k$ let us denote the set of all interacting elements by $\mcT^h_k :=  \{\left.\mcE_\ldx \in \mcT^h \right| (\mcE_k, \mcE_{\ldx}) \in E_S\}$. We then define the subgraph $\mbT_{S_k}:=\left(\mcT^h_k,  E_{S_k} \right)$ of $\mbT_S$ with vertices $\mcT^h_k$ and edges given by
\begin{align*}
    E_{S_k} :=  \left\lbrace \left. (\mcE_k, \mcE_{\ldx}) \in \mcT^h \times \mcT^h  \, \right| \,  (\mcE_k, \mcE_{\ldx}) \in E_S
    \right\rbrace.
\end{align*}

\begin{assumption}
We assume without loss of generality that each of the graphs $\mbT_{S_k}$, $\mbT_{adj}$, and $\mbT_S$ is connected. 
\label{assump:connectedInteraction}
\end{assumption}

\begin{figure}[ht]
    \begin{subfigure}{.5\textwidth}
	    \centering
        \includegraphics[scale=.15]{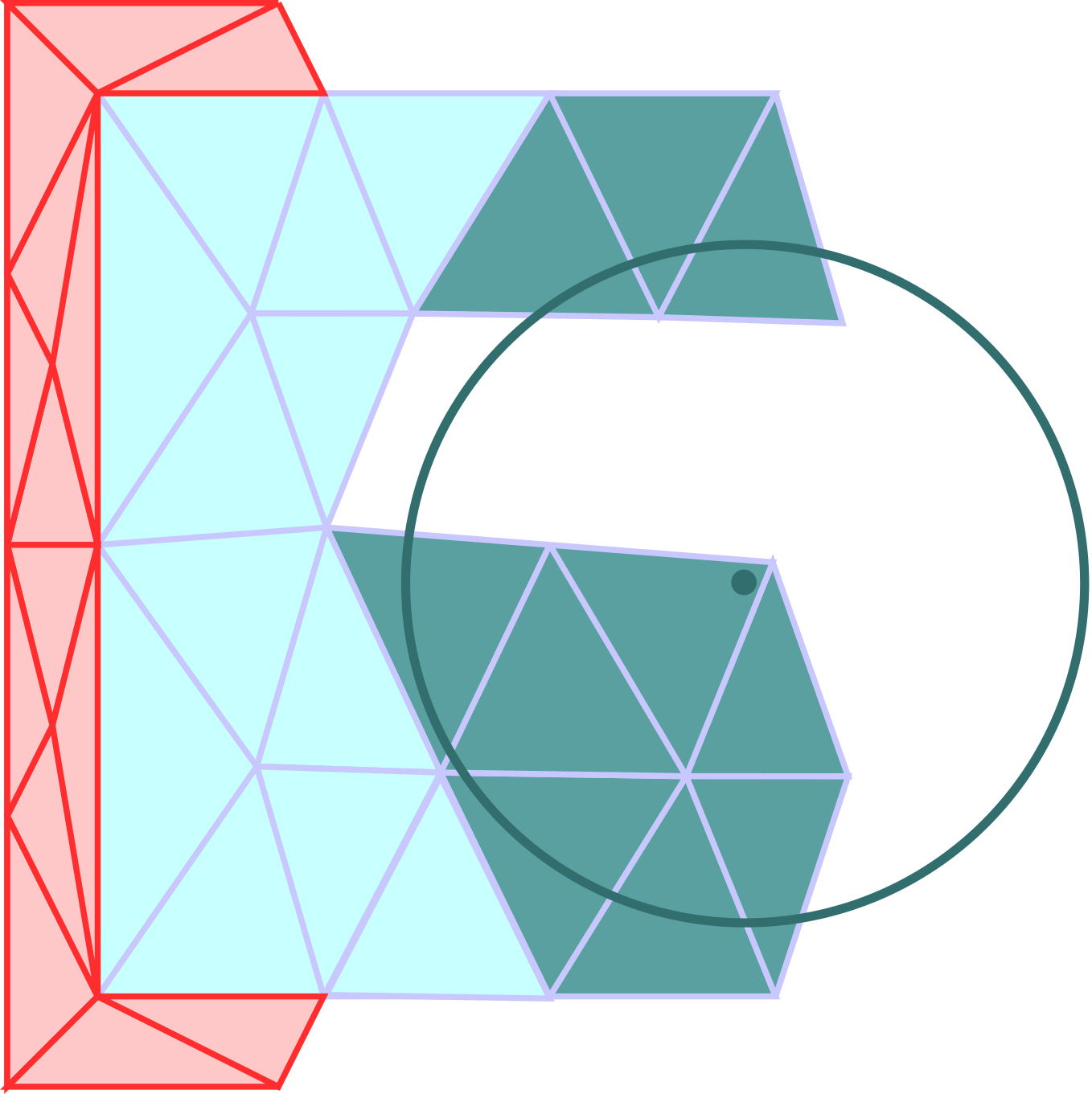}
        \caption{}
        \label{fig:bfsFail}
    \end{subfigure}
    \begin{subfigure}{.49\textwidth}
	    \centering
        \includegraphics[scale=.15]{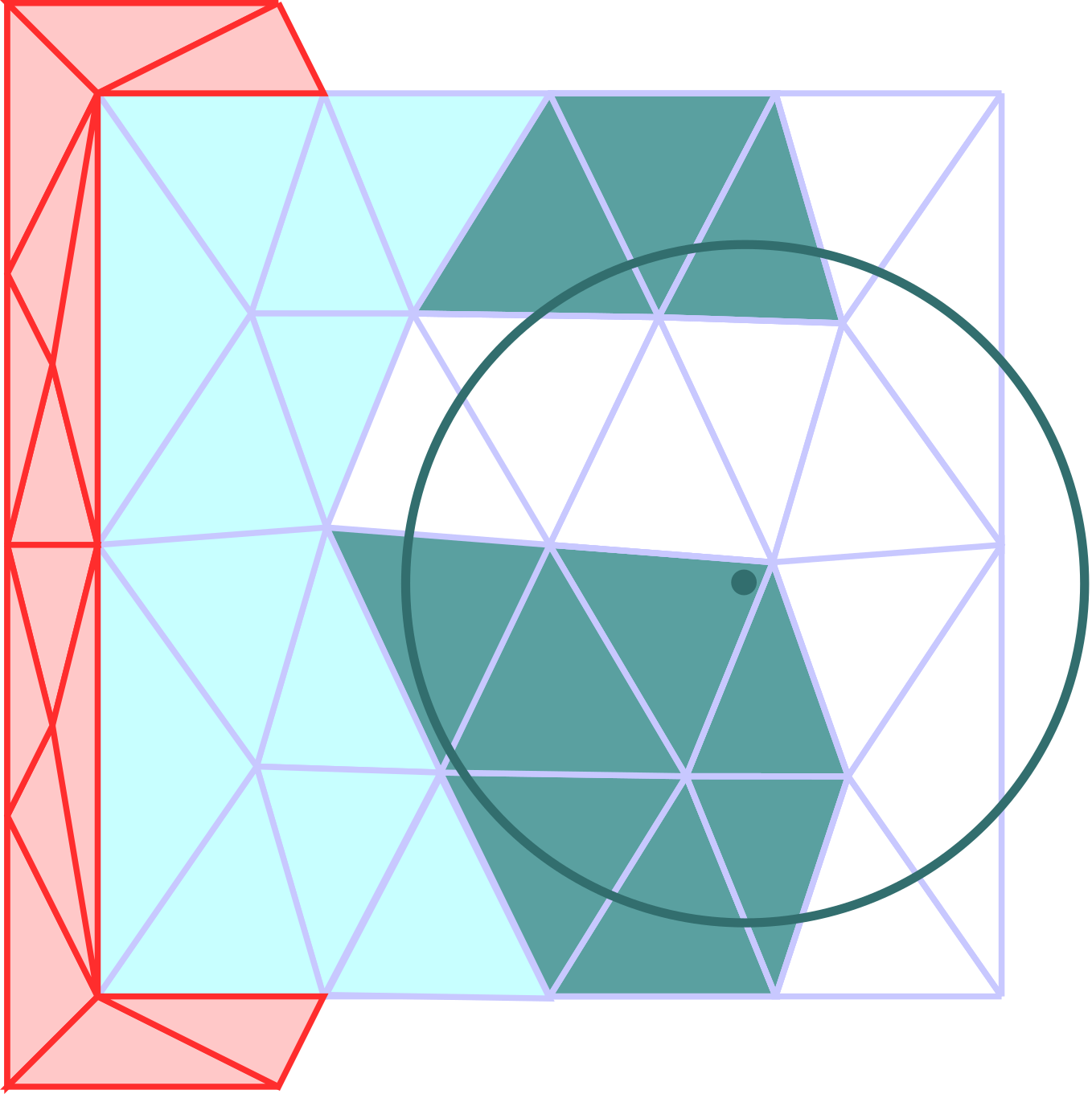}
        \caption{}
        \label{fig:disconnected-domains}
    \end{subfigure}
    \caption{ 
    (a) The interaction neighborhood of an element might be disconnected even if the domain is connected.
    (b) A convex hold all domain $\holdAll$ allows to account for possible interactions
    within disconnected interaction sets.}
\end{figure}

It is clear that Assumption \ref{assump:connectedInteraction} can be violated if $\base$ is not connected or even if it is connected as depicted  in Figure \ref{fig:bfsFail}.
As a remedy, the implementation allows to mark elements which are not part of the finite element discretization.
That way we can add artificial vertices and elements to the mesh. One option then is to connect all vertices at the boundary of the connected components of $\base$ with one artificial vertex, which introduces an adjacency between those elements. 
Another, straightforward option is to embed $\base$ into a bounded and convex \textit{hold-all domain} $\holdAll \subset \R^{\ndim}$, which guarantees that the graphs $\mbT_{S_k}$, $\mbT_{adj}$, and $\mbT_S$ are always connected, see Figure \ref{fig:disconnected-domains}.
In that sense, Assumption \ref{assump:connectedInteraction} does not cause any loss of generality.
If Assumption \ref{assump:connectedInteraction} holds, we can 
recover the subgraph $\mbT_{S_k}$ with a truncated breadth--first traversal of $\mbT_{adj}$ starting
in the root node $\mcE_k$ as given in Algorithm \ref{alg:TraversalofInteratctionSet}.
To that end the root node $\mcE_k$ is appended to a queue $Q := [\mcE_k]$.
While $Q$ is not empty, an element $\mcE_{\widetilde{\ldx}}$ is read and removed  from the queue. Its immediate neighbour, say $\mcE_{\ldx}\in\mcN(\mcE_{\widetilde{\ldx}})$, are obtained from the adjacency graph and the integrals $\widehat{A}_{\kdx\ldx, \idx {\jdx}}^\#$ are successively evaluated for all $\mcE_{\ldx} \in \mcN(\mcE_{\widetilde{\ldx}})$. The element $\mcE_{\ldx}$ is added to the queue whenever the integral does not vanish, and finally marked as visited. This procedure automatically truncates the search to interacting elements for any connected interaction neighborhood.

\begin{algorithm}
	\caption{Traversal of interaction neighborhood} 
	\label{alg:TraversalofInteratctionSet}
	\begin{algorithmic}[1]
		\For {$\mcE_k$ in $\mcT^h$ in parallel threads}
		\State Append $\mcE_k$ to the queue $Q$
		    \While {$Q$ is not empty}
		        \State Read and remove $\mcE_{\widetilde{\ldx}}$ from $Q$ 
				\For {$\mcE_{\ldx}$  in local neighborhood $\mcN (\mcE_{\widetilde{\ldx}})$} 
				    \If {$\mcE_{\ldx}$ is not \textbf{visited}}
				        \State Evaluate  $\localStiffnessMatrix{\kdx}{\ldx}$
				        (Algorithm \ref{alg:populationStiffnessMatrix})
				        \State \textbf{if} $\localStiffnessMatrix{\kdx}{\ldx} \neq 0$ append $\mcE_{\ldx}$ to $Q$
				        \State Mark $\mcE_{\ldx}$ as visited
				    \EndIf
				\EndFor
			\EndWhile
		\EndFor
	\end{algorithmic} 
\end{algorithm}
% NUMERICAL RESULTS
\section{Numerical Examples}\label{sec:numerics}
We solve a truncated fractional--type steady--state diffusion equation, a linear bond--based peridynamics equation \cite{Silling_birthperidym, peridym2} and a steady--state diffusion problem based on the \InfinityNorm ball.
In the first two examples we demonstrate the convergence rate of the approximate solutions to a manufactured solution as the mesh size $h\to 0$. In the latter example we demonstrate the asymptotic compatibility of the discretization scheme for a vanishing horizon $\delta$.

To this end let
$\Omega := (0,0.5)^2 \subset \R^2$. 
We define the interaction domain of $\Omega$ by
$
    \Omega_D := [ -\delta,~ 0.5 + \delta ]^2 ~ \setminus ~ \Omega,
$
so that $\base = [-\delta,0.5+\delta]^2$ is compact.
We then want to solve the nonlocal Dirichlet problem
\begin{align}
\label{prob:nonlocalDirichlet}
    \begin{cases}
    -\mcL_\delta \ub = \fb& \textrm{on } \Omega,\\
    \ub = \gb& \textrm{on } \Omega_D.\\
    \end{cases}
\end{align}
We define the function spaces
\begin{align*}
&V(\base, \R^n) = \lbrace \ub \in L^2(\base, \R^n) : \| \ub \|_V  < \infty \rbrace,\\
&V_c(\base, \R^n) = \lbrace \ub \in V(\base, \R^n) : \ub = 0 \textrm{ in } \Omega_D\rbrace,
\end{align*}
where
\begin{align*}
	\| \ub \|_V^2 = A(\ub, \ub) + \| \ub \|_{L^2(\base)}^2.
\end{align*}
For given data $\fb  \in L^2(\Omega, \R^n)$
and $\gb := \vb_{|\Omega_D}$, where $\vb \in V(\base, \R^n)$, we call $\ub\in V(\base, \R^n)$ the \textit{weak solution} to problem \eqref{prob:nonlocalDirichlet}, if  
\begin{equation}
\label{prob:weakNonlocalDirichlet}
\begin{aligned}
	A(\ub, \vb) &= (\fb, \vb)~~ \textit{for all}~~\vb \in V_c(\base, \R^n),\\
\textit{and}~~ \ub &= \gb \textrm{ in }\Omega_D. 
\end{aligned}
\end{equation}
The wellposedness of problem \eqref{prob:weakNonlocalDirichlet} for various choices of kernel functions can be found, e.g., in
\cite{acta20, du19, nonloc_vector_calc_large, du_analysis_2012, peridym4}.
By exploiting the known values of $\ub$ on the Dirichlet domain, we can rewrite \eqref{prob:weakNonlocalDirichlet} as
\begin{align}
\label{eq:inhomogeneousDirichletData}
    A_{\Omega \Omega}(\ub, \vb) = (\fb, \vb) - A_{\Omega \Omega_D}(\gb, \vb),
\end{align}
where 
\begin{align*}
A_ {\Omega\Omega}(\ub, \vb)
:=&
\int_{\Omega} \int_{\Omega} (\vb(\xb)-\vb(\yb))^\smalltop  \left(\Cb_\delta(\xb, \yb)\ub(\xb)
- \Cb_\delta(\yb, \xb) \ub(\yb) \right)
d\yb  d \xb 
\\
&+
2\int_{\Omega} \int_{\Omega_D} 
\vb(\xb)^\smalltop 
\Cb_\delta(\xb, \yb)\ub(\xb) d\yb  d \xb
\end{align*}
and
\begin{align*}
A_ {\Omega\Omega_D}(\gb, \vb):=  -
2 \int_{\Omega} \int_{\Omega_D} \vb(\xb)^\smalltop  
\Cb_\delta(\yb, \xb) \gb(\yb) 
d\yb  d \xb.
\end{align*} 
In the stiffness matrix the splitting \eqref{eq:inhomogeneousDirichletData} can be naturally obtained by separating the columns corresponding to the degrees of freedom from the columns corresponding to the nodes on the boundary $\Omega_D$. 

\subsection{Truncated Fractional--type Diffusion}
We define the scalar--valued translationally invariant and symmetric kernel function $\gamma^s_\delta: \R^2 \times \R^2 \rightarrow \R$ by
\begin{align}
    \gamma^s_\delta (\xb, \yb) := c_{s,\delta} ~ \frac{1}{|\xb - \yb|^{2+2s}},
\end{align}
where
\begin{align}
    &c_{s, \delta} = \frac{2-2s}{\pi \delta^{2-2s}} ~~\textrm{ and }~~ s \in (0,1). \label{def:fractional-scaling}
\end{align}
Then the scalar--valued {truncated fractional--type diffusion} operator reads as
\begin{align}
    -\mcL^s_\delta u(\xb) := p.v. \int_{B_\delta(\xb)} \gamma^s_\delta (\xb, \yb)  (u(\xb) - u(\yb)) d\yb.
    \label{op:truncatedFractional}
\end{align}
The wellposedness of problem (\ref{prob:weakNonlocalDirichlet}) for this choice of kernel is studied in \cite{SIREV}. 
The constant $c_\delta$ depends on $\delta$ and $s$ and is chosen such that the operator converges to the classical Laplacian as $\delta \to 0$; see, e.g., \cite[Lemma 7.4.1]{Vollmann2019}. For another choice of the constant convergence to the fractional Laplacian as $\delta \to \infty$ can also be obtained \cite{DElia2013TheFL}.

In the example above, we choose the manufactured solution $u(\xb) = \xb_1^2 \xb_2 + \xb_2^ 2$ and set $f(\xb) := -\Delta u(\xb)= -2(\xb_2 + 1)$ in $\Omega$ and $g(\xb):= u(\xb)$ on $\Omega_D$. Since the correctly scaled nonlocal operator equals the classical Laplacian operator on polynomials of order up to three (see, e.g., \cite{Vollmann2019}), we have that $u(\xb)$ is the solution of problem (\ref{prob:weakNonlocalDirichlet}). 
Furthermore, we choose $s:=0.5$, $\delta = 0.2$ and various mesh sizes $h$ as given in the tables below.  
In view of Table \ref{tab:quadratureCases}, for pairs of disjoint elements we use as $\commonQuadRule$ a 7-point quadrature rule\footnote{Specifically, the quadrature points are the barycenter, the vertices, and the mid-side points of the reference triangle $\widehat{\mcE}$ and the corresponding weights are $\frac{27}{60}\!\cdot\!\frac12$, $\frac{3}{60}\!\cdot\!\frac12$, and $\frac{8}{60}\!\cdot\!\frac12$, respectively.} on each, i.e., outer and inner, element. Since six of the seven points are located on the boundary of the triangle, this choice has proven to be advantageous in case of truncated kernels as it mitigates the issue described in Figure \ref{fig:retriangulate-infty}; for more details see \cite{approxBall}. For intersecting pairs we need a quadrature rule on 
$(0,1)^4$ after the integral transformations are performed. For $\duffyQuadRule$ we choose a tensor product of a 5-point Gauss quadrature rule.

The convergence rates on
a continuous Galerkin ansatz space are shown in Table \ref{tab:hto0fractionalCG}.
For our choice of $s = 0.5$ a discontinuous Galerkin ansatz space is also conforming, and the results are presented in Table \ref{tab:hto0fractionalDG}. In both settings, we observe the expected second--order convergence rate as the mesh size $h\to 0$; see, e.g., \cite{SIREV}.
Note that the given examples violate Assumption \ref{assump:connectedInteraction} in the first stage of the experiments as $2 h > \delta$. We see that the first rates in both tables are affected by this.
\begin{table}[ht]
    \centering
    %\subsubsection{Setting}\label{setting}
%
%\begin{longtable}[]{@{}ll@{}}
%\toprule
%Ansatz space & CG\tabularnewline
%Right hand side & linear\tabularnewline
%\textbf{Kernel} & \textbf{fractional}\tabularnewline
%Horizon \(\delta\) & 0.2\tabularnewline
%Fractional constant \(s\)\tabularnewline
%(Default -1) & 0.5\tabularnewline
%\textbf{Intgr. remote pairs} & \textbf{retriangulate}\tabularnewline
%With caps & True\tabularnewline
%Quadrule outer element & 7\tabularnewline
%Quadrule inner element & 7\tabularnewline
%\textbf{Intgr. close pairs}\tabularnewline
%(Relevant only if singular) & \textbf{fractional}\tabularnewline
%Singular quad degree & 5\tabularnewline
%\bottomrule
%\end{longtable}
%
%\subsubsection{Rates}\label{rates}
%
\begin{tabular}{c c c c }%[]{@{}lllll@{}}
\hline
$h$ & dof & L2 Error & Rates \\
\hline
1.41e-01 & 1.60e+01 & 7.19e-04 & 0.00e+00\\
7.07e-02 & 8.10e+01 & 1.65e-04 & 2.13e+00 \\
3.54e-02 & 3.61e+02 & 4.09e-05 & 2.01e+00 \\
1.77e-02 & 1.52e+03 & 1.01e-05 & 2.01e+00 \\
8.84e-03 & 6.24e+03 & 2.41e-06 & 2.07e+00 \\
\hline
\end{tabular}

    \caption{Convergence rates for the truncated fractional diffusion operator (\ref{op:truncatedFractional}), $\delta=0.2$ and $h \to 0$ in a continuous Galerkin ansatz space.}
    \label{tab:hto0fractionalCG}
\end{table}

\begin{table}[ht]
    \centering
    %\subsubsection{Setting}\label{setting}
%
%\begin{longtable}[]{@{}ll@{}}
%\toprule
%Ansatz space & DG\tabularnewline
%Right hand side & linear\tabularnewline
%\textbf{Kernel} & \textbf{fractional}\tabularnewline
%Horizon \(\delta\) & 0.2\tabularnewline
%Fractional constant \(s\)\tabularnewline
%(Default -1) & 0.5\tabularnewline
%\textbf{Intgr. remote pairs} & \textbf{retriangulate}\tabularnewline
%With caps & True\tabularnewline
%Quadrule outer element & 7\tabularnewline
%Quadrule inner element & 7\tabularnewline
%\textbf{Intgr. close pairs}\tabularnewline
%(Relevant only if singular) & \textbf{weakSingular}\tabularnewline
%Singular quad degree & 5\tabularnewline
%\bottomrule
%\end{longtable}
%
%\subsubsection{Rates}\label{rates}
%
\begin{tabular}[]{@{}cccc@{}}
\hline
$h$ & dof & L2 Error & Rates \\
\hline
1.41e-01 & 1.44e+02 & 7.58e-04 & 0.00e+00 \\
7.07e-02 & 5.94e+02 & 1.60e-04 & 2.24e+00 \\
3.54e-02 & 2.39e+03 & 4.06e-05 & 1.98e+00 \\
1.77e-02 & 9.59e+03 & 1.01e-05 & 2.01e+00 \\
8.84e-03 & 3.84e+04 & 2.43e-06 & 2.06e+00 \\
\hline
\end{tabular}

    \caption{Convergence rates for the truncated fractional diffusion operator (\ref{op:truncatedFractional}), $\delta=0.2$ and $h \to 0$ in a discontinuous Galerkin ansatz space.}
    \label{tab:hto0fractionalDG}
\end{table}
\subsection{Bond--Based Peridynamics}
The translationally invariant and symmetric \textit{linear peridynamic} kernel is given by
\begin{align*}
     \Cb_\delta(\xb, \yb) :=  c_{\delta} ~  C(\xb - \yb)~ \ind_{B_\delta(\xb)}(\yb),
\end{align*}
where  
\begin{align}
    C(\xb - \yb) := \frac{(\xb - \yb) (\xb - \yb)^\top}{|\xb - \yb|^3} ~~\text{and}~~ c_\delta := \frac{3}{\delta^3}.
    \label{ker:peridynamics}
\end{align}
The corresponding linear peridynamic operator then reads as
\begin{align}
    -\mcP_\delta \ub(\xb) :=  c_{\delta}\int_{B_\delta(\xb)}C(\xb - \yb)  (\ub(\xb) - \ub(\yb)) d\yb.
    \label{op:peridynamics}
\end{align}
In \cite{peridym4} the wellposedness is established for problem (\ref{prob:weakNonlocalDirichlet}),
where $\fb \in L^2(\Omega, \R^{\ddim})$ and $\gb \in L^2(\Omega_D, \R^{\ddim})$.
For the given constant in \eqref{ker:peridynamics} it is also shown there, that the peridynamic operator $-\mcP_\delta$ converges to the local Navier operator
\begin{align}
    -\mcP_0 \ub (\xb):= -\frac{\pi}{4} \Delta\ub (\xb) - \frac{\pi}{2} \nabla \Div \ub (\xb)
\end{align}
\text{as} $\delta\to 0$. A similar convergence result is also obtained for the corresponding weak solutions. Thus, in the given example we choose the manufactured polynomial solution $\ub(\xb) := (\xb_2^2, \xb_1^2\xb_2)$, set $\fb(\xb) :=-\mcP_0 \ub (\xb)= -\frac{\pi}{2}\left( 1 + 2\xb_1, \xb_2 \right)$ in $\Omega$ and as Dirichlet constraints choose $\gb(\xb) := \ub(\xb)$ on $\Omega_D$. Similarly to the diffusion case above, we again obtain that $\ub(\xb)$ is the solution of \eqref{prob:weakNonlocalDirichlet} due to the correct scaling of the operator \cite{peridym4}. The results are presented in Table \ref{tab:hto0peridynamicsCG} and Table \ref{tab:hto0peridynamicsDG} for a continuous and discontinuous Galerkin ansatz, respectively. In both settings we observe second--order convergence rate as the mesh size $h\to 0$.
\begin{table}[ht]
    \centering
    \begin{tabular}{c c c c c}%[]{@{}lllll@{}}
\hline
$h$ & dof & L2 Error & Rates \\
\hline
1.41e-01 & 3.20e+01 & 7.47e-04 & 0.00e+00 \\
7.07e-02 & 1.62e+02 & 1.82e-04 & 2.04e+00 \\
3.54e-02 & 7.22e+02 & 4.58e-05 & 1.99e+00 \\
1.77e-02 & 3.04e+03 & 1.12e-05 & 2.04e+00 \\
8.84e-03 & 1.25e+04 & 2.63e-06 & 2.09e+00 \\
\hline
\end{tabular}
%rates_kronos_0531_13-51-29

    \caption{Convergence rates for peridynamic operator (\ref{op:peridynamics}), $\delta=0.1$ and $h \to 0$ in a continuous Galerkin ansatz space.}
    \label{tab:hto0peridynamicsCG}
\end{table}

\begin{table}[ht]
    \centering
    \begin{tabular}{c c c c c}%[]{@{}lllll@{}}
\hline
$h$ & dof & L2 Error & Rates \\
\hline
1.41e-01 & 2.88e+02 & 6.64e-04 & 0.00e+00 \\
7.07e-02 & 1.19e+03 & 1.86e-04 & 1.84e+00 \\
3.54e-02 & 4.78e+03 & 4.96e-05 & 1.90e+00 \\
1.77e-02 & 1.92e+04 & 1.25e-05 & 1.99e+00 \\
8.84e-03 & 7.68e+04 & 2.92e-06 & 2.10e+00 \\
\hline
\end{tabular}
%rates_kronos_0531_13-51-29

    \caption{Convergence rates for peridynamic operator (\ref{op:peridynamics}), $\delta=0.1$ and $h \to 0$ in a discontinuous Galerkin ansatz space.}
    \label{tab:hto0peridynamicsDG}
\end{table}

\subsection{Diffusion with Infinity Ball Truncation}
Here, we consider a constant kernel truncated by the \InfinityNorm ball. Specifically, we choose the constant to be $c_\delta^\infty := \frac{3}{4  \delta ^ 4}$, which ensures the convergence to the local Dirichlet problem for vanishing horizon $\delta \to 0$; see, e.g., \cite{Vollmann2019}. The nonlocal operator is then given by
\begin{align}
    -\mathcal{L}^\infty_\delta u (\xb) := 
    c_\delta^\infty
    \int_{B_\delta^\infty(\xb)}  
    (u(\xb) - u(\yb)) d \yb.
    \label{op:infinity}
\end{align}
The truncation by the \InfinityNorm ball is implemented without geometric error, which allows numerical tests of the asymptotic compatibility of the finite element discretization \cite{du_tian_locallim}. 
For the numerical experiment we choose the manufactured solution $u(\xb) = \sin( 4 \pi\xb_1) \sin(4 \pi\xb_2 )$, set 
$f(\xb) := -\Delta u(\xb) = 32 \pi^2 \sin(4 \pi \xb_1) \sin(4 \pi \xb_2 )$ in $\Omega$
and
$g(\xb) := u(\xb)$ on $\Omega_D$. Note that opposed to the previous examples, now for a fixed $\delta >0$ the function $-\mcL_\delta^\infty u$ differs from $-\Delta u$, since $u$ is not chosen to be a polynomial of low degree. Thus the solutions of the nonlocal and the local Dirichlet problem differ from each other and we can better observe the asymptotic compatibility of the discretization scheme.
We run tests for a fixed mesh size $h$ and vanishing $\delta$, see Table \ref{tab:delta-0Linfty}, as well as for a horizon dependent mesh size $h= \sqrt{2}\delta$ and vanishing $\delta$, see Table \ref{tab:delta_h-0Linfty}. In both cases we observe a second--order convergence rate as $\delta \to 0$.
\begin{table}[ht]
    \centering
    %\subsubsection{Setting}\label{setting}
%
%\begin{longtable}[]{@{}ll@{}}
%\toprule
%Ansatz space & CG\tabularnewline
%Right hand side & tensorsin\tabularnewline
%\textbf{Kernel} & \textbf{constantLinf2D}\tabularnewline
%orizon \(\delta\) & 0.0125\tabularnewline
%Fractional constant \(s\)\tabularnewline
%(Default -1) & -1\tabularnewline
%\textbf{Intgr. remote pairs} &
%\textbf{retriangulateLinfty}\tabularnewline
%With caps & True\tabularnewline
%Quadrule outer element & 7\tabularnewline
%Quadrule inner element & 7\tabularnewline
%\textbf{Intgr. close pairs}\tabularnewline
%(Relevant only if singular) &
%\textbf{retriangulateLinfty}\tabularnewline
%Singular quad degree & 1\tabularnewline
%\bottomrule
%\end{longtable}
%
%\subsubsection{Rates}\label{rates}
%
\begin{tabular}{c c c c c}
\hline
$h$ & $\delta$ & L2 Error & Rates \\
\hline
8.84e-03 & 2.00e-01 & 2.03e-01 & 0.00e+00 \\
8.84e-03 & 1.00e-01 & 3.98e-02 & 2.35e+00 \\
8.84e-03 & 5.00e-02 & 8.81e-03 & 2.18e+00 \\
8.84e-03 & 2.50e-02 & 2.08e-03 & 2.08e+00 \\
8.84e-03 & 1.25e-02 & 5.17e-04 & 2.01e+00 \\
\hline
\end{tabular}

    \caption{Convergence rates for the \InfinityNorm ball (\ref{op:infinity}), fixed $h$ and $\delta \to 0$ in a continuous Galerkin ansatz space.}
    \label{tab:delta-0Linfty}
\end{table}

\begin{table}[ht]
    \centering
    %\subsubsection{Setting}\label{setting}
%
%\begin{longtable}[]{@{}ll@{}}
%\toprule
%Ansatz space & CG\tabularnewline
%Right hand side & tensorsin\tabularnewline
%\textbf{Kernel} & \textbf{constantLinf2D}\tabularnewline
%Horizon \(\delta\) & 0.0125\tabularnewline
%Fractional constant \(s\)\tabularnewline
%(Default -1) & -1\tabularnewline
%\textbf{Intgr. remote pairs} &
%\textbf{retriangulateLinfty}\tabularnewline
%With caps & True\tabularnewline
%Quadrule outer element & 7\tabularnewline
%Quadrule inner element & 7\tabularnewline
%\textbf{Intgr. close pairs}\tabularnewline
%(Relevant only if singular) &
%\textbf{retriangulateLinfty}\tabularnewline
%Singular quad degree & 1\tabularnewline
%\bottomrule
%\end{longtable}
%\subsubsection{Rates}\label{rates}
\begin{tabular}{c c c c c}
\hline
$h$ & $\delta$ & L2 Error & Rates \\
\hline
1.41e-01 & 2.00e-01 & 2.00e-01 & 0.00e+00 \\
7.07e-02 & 1.00e-01 & 4.01e-02 & 2.32e+00 \\
3.54e-02 & 5.00e-02 & 8.85e-03 & 2.18e+00 \\
1.77e-02 & 2.50e-02 & 2.10e-03 & 2.07e+00 \\
8.84e-03 & 1.25e-02 & 5.17e-04 & 2.03e+00 \\
\hline
\end{tabular}

    \caption{Convergence rates for the \InfinityNorm ball (\ref{op:infinity}), $\delta = \sqrt{2}h$ and $h,\delta \to 0$ in a continuous Galerkin ansatz space.}
    \label{tab:delta_h-0Linfty}
\end{table}

\subsection{Parallel Complexity of the Assembly process}
The number of elements in each interaction neighborhood grows quadratically in 2d if the diameter of the elements $h$ is decreased for fixed $\delta$. Thus, the assembly of the system matrix with retriangulations as described in the beginning of Section \ref{subsec:quad-truncation} becomes a costly procedure. 
Therefore, a matrix free approach is too expensive and
we store the stiffness matrix in a sparse format.
Furthermore, it makes sense to 
share the work among multiple threads. 
The multithreading is implemented using OpenMP \cite{dagum1998openmp} and the work is shared by a partitioning of finite elements
as given in Algorithm \ref{alg:TraversalofInteratctionSet}, line 1. 
Due to the nonlocality of the operator several threads might need to access identical entries in the global stiffness matrix at the same time to store their contribution. 
OpenMP allows so called critical sections to organize the manipulation of shared variables.
This avoids write conflicts but it would tremendously slow down the computation.
In order to avoid a critical section during the assembly each thread separately allocates its portion of the global stiffness matrix.
The size of the overlap among the submatrices in the threads, i.e., the amount of additional memory requirements due to the parallelization, depends on the nonlocal overlap of the subdomains.
We can therefore reduce the memory requirements by partitioning the domain with metis \cite{karypis1998fast} instead of a scheduler of OpenMP. 
The submatrices are finally added together into a single sparse matrix. 
\begin{table}[ht]
	\centering
	\begin{tabular}{l | ccccccc}
		Threads		        &		1 & 2 & 4 & 8 & 16 & 32 & 64			 \\
		Time [s]		    &	2,153  & 1,169 & 583 & 302 & 166 & 107 & 80 \\
		Parallel efficiency &   -      & 0.92 & 0.92 & 0.89 & 0.81 & 0.63 & 0.42
	\end{tabular}
	\caption{Assembly time for a system with 24,336 degrees of freedom. }
	\label{tab:assemblyScaling}
\end{table}
\begin{figure}[ht]
	\centering
	\includegraphics[scale=.6]{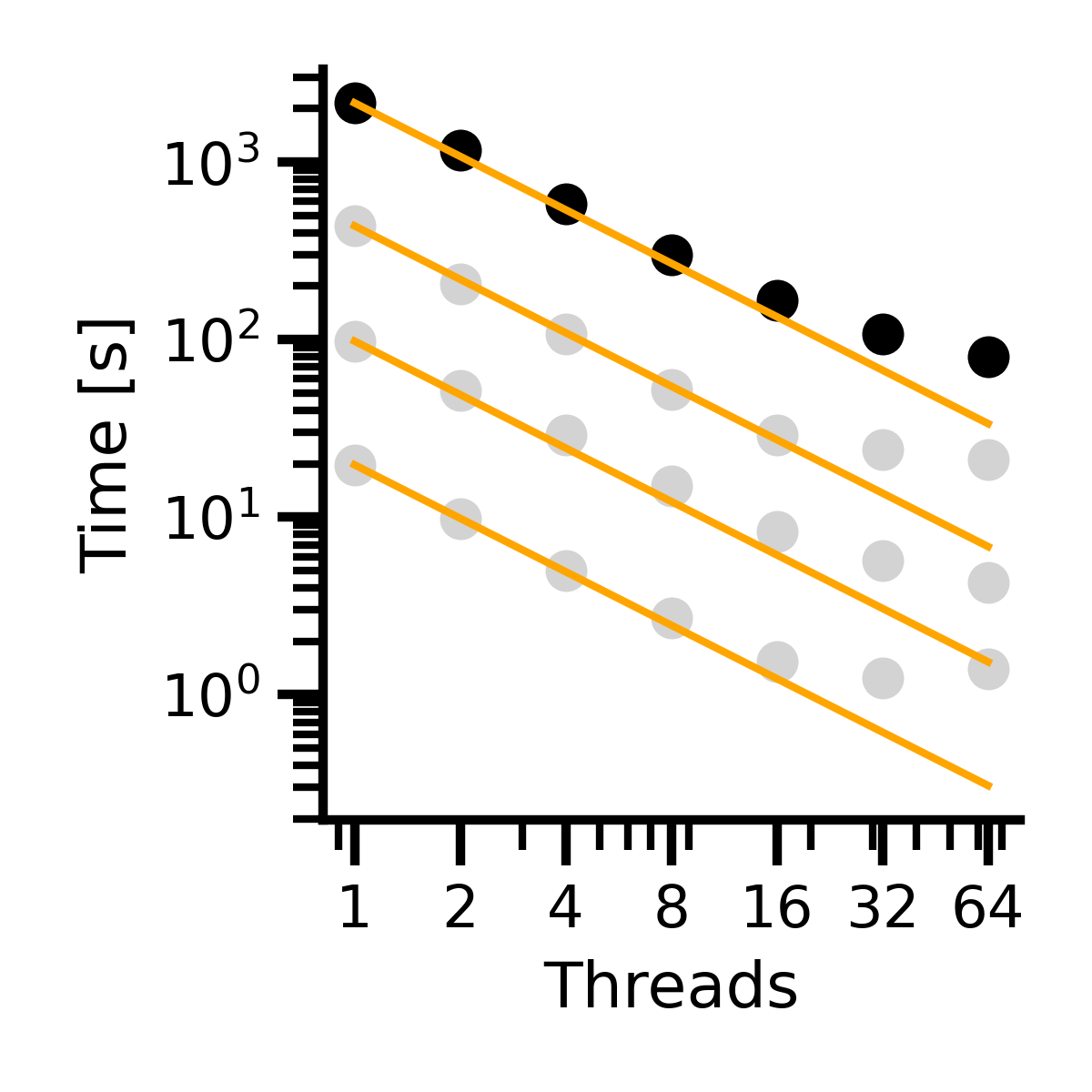}
	\caption{
		Strong scaling for 24,336 (black) and 10,848, 6,120 and 2,736 (gray) degrees of freedom. 
	}
	\label{fig:assemblyScaling}
\end{figure}

We depict the parallel scaling on a computer with two 2.20GHz Intel Xeon CPUs with 22 cores per socket and two threads per core. The machine has 756 GB of RAM.
Table \ref{tab:assemblyScaling} and Figure \ref{fig:assemblyScaling} (black dots) show the run time of an assembly on a regular grid on 
a base domain $\base = [-\delta, 0.5 + \delta]^2$ with mesh size
$h=$ 7.1e-03
and 
$\delta=$ 0.1. The related linear system has 24,336 degrees of freedom and 45,022,167 nonzero entries. The number of threads is increased by a factor up to 64 while the time drops by a factor of 1/27. The scaling looks perfect for up to 16 parallel threads and the effect diminishes from then on. Figure \ref{fig:assemblyScaling} also shows the scaling for smaller problems with 
10,848, 6,120 and 2,736 degrees of freedom (gray dots) which show a similar behavior. 
% CONCLUSION
\section{Conclusion} \label{sec:conclusion}
 The code \texttt{nlfem} is a flexible tool to set up numerical experiments for researchers. The documentation also describes the extension by user defined kernels which allows to consider a large problem class.
It can assemble nonlocal problems in 1d and 3d using the generically implemented barycenter method \cite{approxBall}. 
Furthermore,  \texttt{nlfem} can be used for a recently developed substructuring based domain decomposition method for nonlocal operators \cite{capodaglio2020general} which requires reweighting of the kernel function.
We therefore hope to bridge efficiency and flexibility to obtain a convenient Python package which nourishes the current development in the field of finite element methods for nonlocal operators and enables easy validations of new theory without the effort of implementing code from scratch. \\

%%%%%%%%%%%%%%%%%%%%%%%%%%%%%%%%%%%%%%%%%%%%%%%%%%%%%%%%%%%%%%%%%%%%%%%%%%%%%%%%%
% CLOSING: Acknowledgement, Literature, Appendix
%%%%%%%%%%%%%%%%%%%%%%%%%%%%%%%%%%%%%%%%%%%%%%%%%%%%%%%%%%%%%%%%%%%%%%%%%%%%%%%%%

\pagestyle{plain}
\bibliographystyle{plain}
\bibliography{literature.bib}

\end{document}